\documentclass[10pt]{amsart}  
\usepackage{amsmath,amssymb,amsfonts,graphicx}

\usepackage[colorlinks]{hyperref}

\title{New upper bounds for the density of translative packings of
  three-dimensional convex bodies with tetrahedral symmetry}

\author{Maria Dostert} 
\address{M.~Dostert, Mathematisches Institut, Universit\"at zu
  K\"oln, Weyertal~86--90, 50931 K\"oln, Germany}
\email{m.dostert@uni-koeln.de}

\author{Crist\'obal Guzm\'an}
\address{C.~Guzm\'an, Facultad de Matem\'aticas \& 
Escuela de Ingenier\'ia, Pontificia Universidad Cat\'olica de Chile, 
Av. Vicu\~na Mackenna 4860, Macul, Santiago, Chile}
\email{crguzmanp@uc.cl}

\author{Fernando M\'ario de Oliveira Filho}
\address{F.M.~de Oliveira Filho, 
Instituto de Matem\'atica e Estat\'\i stica, 
Universidade de S\~ao Paulo, Rua do Mat\~ao, 1010, 05508-090 S\~ao
Paulo, Brazil}
\email{fmario@gmail.com}

\author{Frank Vallentin} 
\address{F.~Vallentin, Mathematisches Institut, Universit\"at zu
  K\"oln, Weyertal~86--90, 50931 K\"oln, Germany}
\email{frank.vallentin@uni-koeln.de}

\date{January 10, 2017}

\thanks{The fourth author was partially supported by VIDI grant
  639.032.917 from the Netherlands Organization for Scientific
  Research (NWO)}



\newcommand{\defi}[1]{\textit{#1}}

\newcommand{\R}{\mathbb{R}}
\newcommand{\N}{\mathbb{N}}

\newcommand{\C}{\mathbb{C}}

\newtheorem{defin}{Definition}[section]

\newtheorem{proposition}[defin]{Proposition}
\newtheorem{theorem}[defin]{Theorem}

\DeclareMathOperator{\ort}{O}
\DeclareMathOperator{\vol}{vol}

\DeclareMathOperator{\Hom}{Hom}
\DeclareMathOperator{\Harm}{Harm}

\DeclareMathOperator{\GL}{GL}

\newcommand{\Acal}{\mathcal{A}}
\newcommand{\Bcal}{\mathcal{B}}
\newcommand{\Ccal}{\mathcal{C}}
\newcommand{\Dcal}{\mathcal{D}}
\newcommand{\Fcal}{\mathcal{F}}

\newcommand{\Ical}{\mathcal{I}}
\newcommand{\Scal}{\mathcal{S}}
\newcommand{\Kcal}{\mathcal{K}}

\newcommand{\bt}{{\mathsf{B}_3}}
\newcommand{\btd}{{\widehat{\mathsf{B}}_3}}

\newcommand{\tp}{{\sf T}}

\begin{document}

\begin{abstract} 
  In this paper we determine new upper bounds for the maximal density
  of translative packings of superballs in three dimensions (unit
  balls for the $l^p_3$-norm) and of Platonic and Archimedean
  solids having tetrahedral symmetry. Thereby, we improve Zong's recent upper bound for the maximal
  density of translative packings of regular tetrahedra from
  $0.3840\ldots$ to $0.3745\ldots$, getting closer to the best known
  lower bound of $0.3673\ldots$

  We apply the linear programming bound of Cohn and Elkies which
  originally was designed for the classical problem of densest packings of
  round spheres. The proofs of our new upper bounds are computational
  and rigorous. Our main technical contribution is the use of
  invariant theory of pseudo-reflection groups in polynomial
  optimization.
\end{abstract}

\maketitle

\markboth{M. Dostert, C. Guzm\'an, F.M.~de Oliveira Filho,
  F.~Vallentin}{New upper bounds for the density of translative
  packings of three-dimensional convex bodies}

\vspace*{-1cm}

\tableofcontents

\section{Introduction}
\label{sec:introduction}

The most famous geometric packing problem is Kepler's conjecture from
1611: The density of any packing of equal-sized spheres into
three-dimensional Euclidean space is never greater than
$\pi/\sqrt{18} = 0.7404 \ldots$ This density is achieved for example
by the ``cannonball'' packing.  In 1998 Hales and Ferguson solved
Kepler's conjecture. Their proof is extremely complicated, involving
more than 200 pages, intensive computer calculations, and checking of
more than 5,000 subproblems. They wrote a book~\cite{Hales2011a} that
contains the entire proof together with supporting material and
commentary.

Very little is known if one goes beyond three-dimensional packings of spheres to packings
of nonspherical objects. Considering nonspherical objects is
interesting for many reasons. For example, using nonspherical objects
one can model physical granular materials accurately. On the other
hand, the mathematical difficulty increases substantially when one
deals with nonspherical objects.

Jiao, Stillinger, and Torquato \cite{Jiao2009a} consider packings of
three-dimensional superballs, which are unit balls of the
$l^p_3$-norm, with $p \geq 1$:
\[
B^p_3 = \{(x_1,x_2,x_3) \in \R^3 : |x_1|^p + |x_2|^p + |x_3|^p \leq 1\}.
\]
Three-dimensional superballs can be synthesized experimentally as
colloids, see Rossi et al.\ \cite{Rossi2011a}. The name ``superball''
is attributed to the Danish inventor Piet Hein who used a
``superellipse'' with $p = 5/2$ in a design challenge of the
redevelopment of the public square Sergels Torg in Stockholm (see also
Rush and Sloane \cite{Rush1987a}, and Gardner \cite{Gardner1989a}).
Life Magazine \cite{Hicks1966a} quotes Piet Hein:
\begin{quote}
\small
  Man is the animal that draws lines which he himself then stumbles
  over. In the whole pattern of civilization there have been two
  tendencies, one toward straight lines and rectangular patterns and
  one toward circular lines. There are reasons, mechanical and
  psychological, for both tendencies. Things made with straight lines
  fit well together and save space. And we can move
  easily---physically or mentally---around things made with round
  lines. But we are in a straitjacket, having to accept one or the
  other, when often some intermediate form would be better.
\end{quote}

%, Life Magazine, Vol. 61 No. 16, 10/14/66, pp. 55–66
% http://books.google.de/books?id=lFYEAAAAMBAJ&lpg=PP1&hl=de&pg=PA55#v=onepage&q&f=false

Back to the work of Jiao, Stillinger, and Torquato
\cite{Jiao2009a}. They construct the densest known packings of $B^p_3$
for many values of $p$. As motivation for their study Jiao, Stillinger,
and Torquato write:
\begin{quote}
\small
  Understanding the organizing principles that lead to the densest
  packings of nonspherical particles that do not tile space is of
  great practical and fundamental interest. Clearly, the effect of
  asphericity is an important feature to include on the way to
  characterizing more fully real dense granular media.

\smallskip

[...]

\smallskip

On the theoretical side, no results exist that rigorously prove the
densest packings of other congruent non-space-tiling particles in
three dimension.
\end{quote}

Torquato and Jiao \cite{Torquato2009a, Torquato2009b} extend the work
on superballs to nonspherical non-differentiable shapes.  They found
dense packings of Platonic and of Archimedean solids.

\begin{quote}
\small
Very little is known about the densest packings of polyhedral
particles that do not tile space, including the majority of the
Platonic and Archimedean solids studied by the ancient Greeks. The
difficulty in obtaining dense packings of polyhedra is related to
their complex rotational degrees of freedom and to the non-smooth
nature of their shapes.
\end{quote}

The optimal, densest lattice packing of each Platonic or Archimedean
solid is known. Minkowski \cite{Minkowski1904a} determines the densest
lattice packing of regular octahedra. Hoylman \cite{Hoylman1970a} uses
Minkowski's method to determine the densest lattice packing of regular
tetrahedra. Betke and Henk \cite{Betke2000a} turn Minkowski's method
into an implementable algorithm and find the densest lattice packings
of all remaining Platonic and all Archimedean solids. Only two of the
Platonic and Archimedean solids are not centrally symmetric, namely
the tetrahedron and the truncated tetrahedron. These are also the only
cases where Torquato and Jiao could use the extra freedom of rotating
the solids to find new packings which are denser than the
corresponding densest lattice packings. Also the dense superball
packings of Jiao, Stillinger, and Torquato are lattice packings. Based
on this evidence they formulate the following conjecture:

\begin{quote}
\small
The densest packings of the centrally symmetric Platonic and
Archimedean solids are given by their corresponding optimal lattice
packings. This is the analogue of Kepler's sphere conjecture for these
solids.
\end{quote}

For a convex body $\mathcal{K}$ in $\mathbb{R}^n$ it is natural to
consider three kinds of increasingly restrictive packings: congruent
packings, translative packings, and lattice packings. A packing of
congruent copies of $\mathcal{K}$ has the form
\[
\mathcal{P} = \bigcup_{i \in \mathbb{N}} (x_i + A_i \mathcal{K}), \quad
\text{with } (x_i, A_i) \in \mathbb{R}^n \times \mathrm{SO}(n), \; i
\in \mathbb{N},
\]
where
$x_i + A_i \mathcal{K}^\circ \cap x_j + A_j \mathcal{K}^\circ =
\emptyset$
whenever the indices $i$ and $j$ are distinct. Here,
$\mathcal{K}^\circ$ denotes the topological interior of the
body~$\mathcal{K}$ and
\[
\mathrm{SO}(n) = \{A \in \mathbb{R}^{n \times n} : AA^{\sf T} = I_n,
\, \det A = 1\}
\]
denotes the special orthogonal group, an index-$2$ subgroup of the
orthogonal group
\[
\mathrm{O}(n) = \{A \in \mathbb{R}^{n \times n} : AA^{\sf T}  = I_n\}.
\]
The (upper) density of $\mathcal{P}$ is
\[
\delta(\mathcal{P}) = \limsup_{r \to \infty} \sup_{c \in \R^n} \frac{\vol (B(c, r) \cap \mathcal{P})}{\vol
  B(c, r)},
\]
where~$B(c, r)$ is the Euclidean ball of radius~$r$ centered at~$c$.
A packing $\mathcal{P}$ is called a translative packing if each
rotation $A_i$ is identity.  A translative packing is called a
lattice packing if the set of $x_i$'s forms a lattice.

Lattice packings are restrictive and many results about them are
known. This is not the case for translative and congruent
packings. While the conjecture of Torquato and Jiao ultimately aims at
congruent packings, the objective of our paper is to develop tools
coming from mathematical optimization which will be useful to make
progress on the conjecture restricted to translative packings. In
particular we prove new upper bounds for the density of densest
translative packings of three-dimensional superballs and of Platonic
and Archimedean solids with tetrahedral symmetry. We use the following
theorem of Cohn and Elkies \cite{Cohn2003a} for this.

\begin{theorem}
\label{thm:cohn-elkies}
Let $\Kcal$ be a convex body in $\R^n$ and let $f \colon \R^n \to \R$
be a continuous $L^1$-function (a continuous function whose absolute
value is Lebesgue integrable). Define the Fourier transform of
$f$ at $u$ by
\[
\widehat{f}(u) = \int_{\R^n} f(x) e^{-2\pi i u \cdot x}\, dx.
\]
Suppose $f$ satisfies the following conditions
\begin{enumerate}
\item[(i)] $\widehat{f}(0) \geq 1$,
\item[(ii)] $f$ is of positive type, i.e.\ $\widehat{f}(u) \geq 0$ for
  every $u \in \R^n$,
\item[(iii)] $f(x) \leq 0$ whenever $\Kcal^\circ \cap (x +
  \Kcal^\circ) = \emptyset$.
\end{enumerate}
Then the density of any packing of translates of $\Kcal$ in~$\R^n$ is
at most $f(0) \vol\Kcal$. 
\end{theorem}

One can find a proof of this theorem in Cohn and Kumar
\cite{Cohn2007a} or for the more general case of translative packings
of multiple convex bodies $\Kcal_1, \ldots, \Kcal_N$ in de Laat,
Oliveira, and Vallentin \cite{Laat2014a}. Originally, Cohn and Elkies
\cite{Cohn2003a} state the theorem only for admissible functions;
these are functions for which the Poisson summation formula applies.

The Cohn-Elkies bound provides the basic framework for proving the
best known upper bounds for the maximum density of sphere packings in
dimensions~$4$, \dots,~$36$.  For some time it was conjectured to
provide tight bounds in dimensions~$8$ and~$24$ and there was very
strong numerical evidence to support this conjecture, see Cohn and
Miller~\cite{Cohn2016a}. However, the only thing missing was a
rigorous proof. Recently, in March 2016, such a proof was found by
Viazovska \cite{Viazovska2016a} for dimension $8$ and a few days
later, building on Viazovska's breakthrough result, by Cohn, Kumar,
Miller, Radchenko, and Viazovska \cite{Cohn2016b} for dimension
$24$. Here the explicit construction of optimal functions $f$ uses the
theory of quasimodular forms from analytic number theory.

De Laat, Oliveira, and Vallentin~\cite{Laat2014a} have proposed a
strengthening of the Cohn-Elkies bound and computed better upper
bounds for the maximum density of sphere packings in dimensions~$4$,
$5$, $6$, $7$, and~$9$.

In all these calculations one can restrict the function $f$ to be a
radial function (a function whose values $f(x)$ only depend on the
norm $\|x\|$ of the vector $x \in \mathbb{R}^n$) because of the
rotational symmetry of the sphere. For the case of packings of
nonspherical objects Cohn and Elkies \cite{Cohn2003a} write:

\begin{quote}
  \small Unfortunately, when [the body we want to pack] is not a
  sphere, there does not seem to be a good analogue of the reduction
  to radial functions in Theorem [\ref{thm:cohn-elkies}]. That makes
    these cases somewhat less convenient to deal with.
\end{quote}

Until now, the Cohn-Elkies bound has only been computed for packings
of spheres. In this paper we show how to apply the Cohn-Elkies bound
for nonspherical objects.

\subsection{New upper bounds for translative packings}

Before we describe our methods we report on the new upper bounds we
obtained and compare them to the known lower and upper bounds. We give
the new upper bounds for three-dimensional superball packings in
Table~\ref{table:superballs}, the new upper bounds for Platonic and
Archimedean solid with tetrahedral symmetry are in
Table~\ref{table:polytopes}.

\renewcommand{\thetable}{\arabic{table}}

\begin{table}[htb]

\begin{center}
\begin{tabular}{lcc}
Body & \multicolumn{2}{c}{Lattice packing} \\
& lower bound & upper bound\\
\hline
$B_3^1$ & 
$18/19 = 0.9473\dots$ \cite{Minkowski1904a}&
 $18/19$ \cite{Minkowski1904a}\\
$B_3^2$ &
$\pi/\sqrt{18} = 0.7404\dots$&     
$\pi/\sqrt{18}$ \cite{Gauss1840a}
\\
$B_3^3$ &   
$0.8095\dots$ \cite{Jiao2009a} & 
$\mathit{0.8236\dots}$\\
$B_3^4$ &  
$0.8698\dots$\cite{Jiao2009a} & 
$\mathit{0.8742\dots}$\\
$B_3^5$ &    
$0.9080\dots$ \cite{Jiao2009a}& 
$\mathit{0.9224\dots}$\\
$B_3^6$ &   
$0.9318\dots$ \cite{Jiao2009a}& 
$\mathit{0.9338\dots}$\\
\hline
\end{tabular}
\end{center}

\bigskip

\begin{center}
\begin{tabular}{lcc}
& \multicolumn{2}{c}{Translative packing}\\
& lower bound & upper bound\\
\hline
$B_3^1$ & 
 $18/19 = 0.9473\dots$ \cite{Minkowski1904a}&
 $\mathit{0.9729\dots}$\\
$B_2^2$ &
$\pi/\sqrt{18} = 0.7404\dots$ & 
$\pi/\sqrt{18}$ \cite{Hales2011a}\\
$B_3^3$ &   
$0.8095\dots$ \cite{Jiao2009a}& 
$\mathit{0.8236\dots}$\\
$B_3^4$ &  
$0.8698\dots$ \cite{Jiao2009a}& 
$\mathit{0.8742\dots}$\\
$B_3^5$ &    
$0.9080\dots$ \cite{Jiao2009a}& 
$\mathit{0.9224\dots}$\\
$B_3^6$ &   
$0.9318\dots$ \cite{Jiao2009a}& 
$\mathit{0.9338\dots}$\\
\hline
\end{tabular}
\end{center}

\bigskip

\begin{center}
\begin{tabular}{lcccccc}
& \multicolumn{2}{c}{Congruent packing}\\
& lower bound & upper bound\\
\hline
$B_3^1$ & 
 $18/19 = 0.9473\dots$ \cite{Minkowski1904a} &
$1 - 1.4\ldots\cdot 10^{-12}$ \cite{Gravel2011a} \\
$B_2^2$ &
$\pi/\sqrt{18} = 0.7404\dots$ & 
$\pi/\sqrt{18}$ \cite{Hales2011a} \\
$B_3^3$ &   
$0.8095\dots$ \cite{Jiao2009a}& 
$<1$\\
$B_3^4$ &  
$0.8698\dots$ \cite{Jiao2009a}& 
$<1$\\
$B_3^5$ &    
$0.9080\dots$ \cite{Jiao2009a}& 
$<1$ \\
$B_3^6$ &   
$0.9318\dots$ \cite{Jiao2009a}& 
$<1$ \\
\hline
\end{tabular}
\end{center}
\bigskip

\caption{Best known bounds for packings of three-dimensional
  superballs. New bounds obtained in this paper are written in italics.}
\label{table:superballs}
\end{table}

\begin{table}[htb]

\begin{center}
\begin{tabular}{lcc}
Body & \multicolumn{2}{c}{Lattice packing}\\
& lower bound & upper bound\\
\hline
Tetrahedron & 
$18/49 =  0.3673\dots$ \cite{Groemer1962a} &       
$18/49$ \cite{Hoylman1970a}\\
Truncated tetrahedron & 
$0.6809\dots$ \cite{Betke2000a} &
$0.6809\dots$ \cite{Betke2000a}\\
Truncated cuboctahedron &  
$0.8493\dots$ \cite{Betke2000a} &
$0.8493\dots$ \cite{Betke2000a}\\
Rhombicuboctahedron &  
$0.8758\dots$ \cite{Betke2000a} & 
$0.8758\dots$ \cite{Betke2000a}\\
Cuboctahedron &  
$0.9183\dots$ \cite{Groemer1962a} &
$0.9183\dots$ \cite{Hoylman1970a}\\
Truncated cube &  
$0.9737\dots$ \cite{Betke2000a} &
$0.9737\dots$ \cite{Betke2000a}\\
\hline
\end{tabular}
\end{center}

\bigskip

\begin{center}
\begin{tabular}{lcc}
& \multicolumn{2}{c}{Translative packing}\\
& lower bound & upper bound\\
\hline
Tetrahedron & 
$18/49 = 0.3673\dots$ \cite{Groemer1962a} &
$\mathit{0.3745\dots}$\\
Truncated tetrahedron & 
$0.6809\dots$ \cite{Betke2000a} & 
$\mathit{0.7292\dots}$ \\
Truncated cuboctahedron &  
$0.8493\dots$ \cite{Betke2000a}&
$0.8758\dots$ \cite{Torquato2009b} \\
Rhombicuboctahedron &  
$0.8758\dots$ \cite{Betke2000a}&
$0.8758\dots$ \cite{deGraaf2011a} \\
Cuboctahedron &  
$0.9183\dots$ \cite{Groemer1962a}&
$\mathit{0.9364\dots}$\\
Truncated cube &  
$0.9737\dots$ \cite{Betke2000a}&
$\mathit{0.9845\dots}$\\
\hline
\end{tabular}
\end{center}

\bigskip

\begin{center}
\begin{tabular}{lcccccc}
& \multicolumn{2}{c}{Congruent packing}\\
& lower bound & upper bound\\
\hline
Tetrahedron & 
$4000/4671 = 0.8563\dots$ \cite{Chen2010a} &
$1 - 2.6\ldots\cdot 10^{-25}$ \cite{Gravel2011a} \\
Truncated tetrahedron & 
$207/208 = 0.9951\dots$ \cite{Jiao2011b}, \cite{Damasceno2012a}& 
$<1$ \\
Truncated cuboctahedron &  
$0.8493\dots$ \cite{Betke2000a} &
$0.8758\dots$ \cite{Torquato2009b} \\
Rhombicuboctahedron &  
$0.8758\dots$ \cite{Betke2000a} & 
$0.8758\dots$ \cite{deGraaf2011a} \\
Cuboctahedron &  
$0.9183\dots$ \cite{Groemer1962a} & 
$<1$\\
Truncated cube &  
$0.9737\dots$ \cite{Betke2000a} &
$<1$\\
\hline
\end{tabular}
\end{center}

\bigskip

\caption{Best known bounds for packings of three-dimensional
  Platonic and Archimedean solids with tetrahedral symmetry.  The
  octahedron, the cube and the truncated octahedron are omitted. New bounds obtained in this paper are written in italics.}
\label{table:polytopes}
\end{table}

\subsubsection{Three-dimensional superballs}

Jiao, Stillinger, and Torquato \cite{Jiao2009a, Jiao2011a}
find dense packings of superballs $B^p_3$ for all values of
$p \geq 1$. Although they principally allow congruent packings in
their computer simulations, the dense packings they find are all
lattice packings. They subdivide the range $p \in [1,\infty)$ into four
different regimes
\[
p \in [1, 2 \ln 3 / \ln 4 = 1.5849\ldots] \cup [1.5849\ldots, 2] \cup [2, 2.3018\ldots] \cup [2.3018\ldots,\infty)
\]
and give for each regime a family of lattices determining dense
packings.  When $p = 1$ then $B^p_3$ is simply the regular octahedron,
a Platonic solid. The optimal lattice packing of regular octahedra has
been determined by Minkowski~\cite{Minkowski1904a}.

Recently, Ni, Gantapara, de Graaf, van Roij, and
Dijkstra~\cite{Ni2012a} claimed that for values of $p$ lying in the
interior of the first and respectively of the second regime, the
packings of Jiao, Stillinger, and Torquato can be improved.

When $p = 2$, then $B^p_3$ is the round unit ball of Euclidean space.
The optimal lattice packing of $B^2_3$ has been determined by
Gauss~\cite{Gauss1840a} using reduction theory of positive quadratic
forms. Here, translative and congruent packings coincide because of
the rotational symmetry of $B^2_3$. Hales~\cite{Hales2011a} proved the
optimality of the cannonball packing among all congruent packings. One
should note that there is an uncountable family of non-lattice
packings achieving the same density. The best upper bound obtainable from
Theorem~\ref{thm:cohn-elkies} is $0.7797\dots$

For $p \geq 2.3018\ldots$ the densest known superball packings are given by
the family of $C_1$-lattices which is defined by
$\mathbb{Z} b_1 + \mathbb{Z} b_2 + \mathbb{Z} b_3$ with
\[
b_1 = (2^{1-\frac{1}{p}}, 2^{1-\frac{1}{p}}, 0),\; b_2 = (2^{1-\frac{1}{p}}, 0,
2^{1-\frac{1}{p}}),\; b_3 = (2s+2^{1-\frac{1}{p}}, -2s, -2s),
\]
where $s$ is the smallest positive root of the equation
\[
(s + 2^{-\frac{1}{p}})^{p} + 2s^{p} -1 = 0,
\]
having density
\[
\frac{\vol B^p_3}{2^{3-\frac{2}{p}}(3s+2^{-\frac{1}{p}})}.
\]

\subsubsection{Platonic and Archimedean solids with tetrahedral symmetry}

We prove a new upper bound for the density of densest translative
packings of regular tetrahedra and improve the upper bound of
$0.3840\ldots$ recently obtained by Zong \cite{Zong2014a} to
$0.3745\ldots$ Groemer~\cite{Groemer1962a} shows that there is a
lattice packing of regular tetrahedra which has density $0.3673\ldots$
and Hoylman~\cite{Hoylman1970a} proves the optimality of Groemer's
packing among lattice packings using Minkowski's method. Finding dense
congruent packings of regular tetrahedra is fascinating. In fact, it
is part of Hilbert's 18th problem. We refer to Lagarias and Zong
\cite{Lagarias2012a} for the history of the tetrahedra packing problem
and to Ziegler~\cite{Ziegler2011a} for an overview on a race for the
best construction.

As a corollary of our new bound for the tetrahedron, we improve Zong's
bound for densest translative packings of the cuboctahedron from
$0.9601\ldots$ to $0.9364\ldots$ This follows from Minkowski's
observation that
\[
\bigcup_{i \in \mathbb{N}} (x_i + \mathcal{K})
\]
is a translative packing of $\mathcal{K}$ if and only if
\[
\bigcup_{i \in \mathbb{N}} \left(x_i + \frac{1}{2}(\mathcal{K}
-\mathcal{K})\right)
\]
is a packing of $\frac{1}{2}(\mathcal{K} - \mathcal{K})$, with
\[
\mathcal{K} - \mathcal{K} = \{x - y : x,y \in \mathcal{K}\}
\]
denoting the Minkowski difference of the body $\mathcal{K}$ with
itself.  The Minkowski difference of a regular tetrahedron with itself
is the cuboctahedron whose volume is $2^3 \cdot 5/2$ times the volume
of the regular tetrahedron.

We omitted the octahedron from Table~\ref{table:polytopes} because it
is $B^1_3$ in Table~\ref{table:superballs}. We also omitted the
cube and the truncated octahedron because with both solids one
can tile three-dimensional space.

\subsection{Computational strategy}

In this section we give a high level description of how we found
suitable functions $f$ for Theorem~\ref{thm:cohn-elkies} for proving
the new upper bounds.

The symmetry group of a convex body $\Kcal \subseteq \R^n$ is
\[
S(\Kcal) = \{A \in \ort(n) : A\Kcal = \Kcal\},
\]
and when considering functions $f$ for Theorem~\ref{thm:cohn-elkies}
the symmetry group of the Minkowski difference $\Kcal - \Kcal$ will be
useful. For $A \in S(\Kcal - \Kcal)$ we have
\[
\begin{split}
& A^{-1} x + \Kcal^\circ \cap \Kcal^\circ = \emptyset 
 \Longleftrightarrow
x + A\Kcal^\circ \cap A\Kcal^\circ = \emptyset \Longleftrightarrow x
\not\in A\Kcal^\circ - A\Kcal^\circ \\
& \quad \Longleftrightarrow x
\not\in A(\Kcal^\circ - \Kcal^\circ)
\Longleftrightarrow x \not\in \Kcal^\circ - \Kcal^\circ
\Longleftrightarrow x + \Kcal^\circ \cap \Kcal^\circ = \emptyset.
\end{split}
\]
Hence we may assume without loss of generality that the function $f$
we are seeking is invariant under the left action of
$S(\Kcal - \Kcal)$, i.e.\
\[
f(A^{-1} x) =  f(x) \quad \text{for all } A \in S(\Kcal - \Kcal).
\]
This assumption reduces the search space and also makes the third
constraint $f(x) \leq 0$ whenever
$x + \Kcal^\circ \cap \Kcal^\circ = \emptyset$ easier to model.

In the case of $\Kcal$ being a superball $B^p_3$, with $p \geq 1$ and
$p \neq 2$, the Minkowski difference $\Kcal - \Kcal$ is $2 B^p_3$ and
its symmetry group is a finite subgroup of the orthogonal group. It is
the octahedral group (which is the same as the symmetry group of the
regular cube $[-1,+1]^3$), which has $48$ elements. The octahedral
group is the reflection group $\mathsf{B}_3$ which is generated by the
three matrices
\begin{equation}
\label{eq:b3-generators}
\begin{pmatrix}
-1 & 0 & 0\\
0 & 1 & 0\\
0 & 0 & 1
\end{pmatrix},\;
\begin{pmatrix}
1 & 0 & 0\\
0 & 0 & -1\\
0 & -1 & 0
\end{pmatrix},\; 
\begin{pmatrix}
0 & 0 & 1\\
0 & 1 & 0\\
1 & 0 & 0
\end{pmatrix},
\end{equation}
where the first one is the reflection at the plane $x_1 = 0$, the
second one is the reflection at the plane $x_2 + x_3 = 0$, and the
last one is the reflection at the plane $x_1 - x_3 = 0$.

In the case of $\Kcal$ being a Platonic or Archimedean solid with
tetrahedral symmetry, the symmetry group of the Minkowski difference
$\Kcal - \Kcal$ is the octahedral group, too.

We specify the function~$f\colon \R^3 \to \R$ via its Fourier
transform~$\widehat{f}$.  If $f$ is invariant under the action of
$S(\Kcal - \Kcal)$ then the same is true for its Fourier transform
$\widehat{f}$. Let~$g$ be a polynomial. We use the following template
for the Fourier transform of~$f$:
\begin{equation}
\label{eq:fhat}
\widehat{f}(u) = g(u) e^{-\pi \|u\|^2}.
\end{equation}
So~$\widehat{f}$ is a Schwartz function (i.e.\ all derivates
$D^\beta \widehat{f}(x)$ exist for all $x \in \mathbb{R}^n$ and all
$\beta \in \mathbb{N}^n$ and
$\sup\{|x^\alpha D^\beta \widehat{f}(x)| : x \in \mathbb{R}^n\} < \infty$ holds for all
$\alpha, \beta \in \mathbb{N}^n$), implying that also~$f$ is a
Schwartz function. In particular, $f$ will be a continuous
$L^1$-function.

If $f$ is invariant under $\mathsf{B}_3$ then so is the polynomial $g$
which specifies $\widehat{f}$. This means that $g$ lies in the ring of
invariants of the group $\mathsf{B}_3$ which is, by the theory of
finite reflection groups, known to be freely generated by three basic
invariants
\begin{equation}
\label{eq:basic-invariants}
\theta_1 = x_1^2 + x_2^2 + x_3^2, \;\; \theta_2 = x_1^4 + x_2^4 + x_3^4, \;\; \theta_3 = x_1^6 + x_2^6 + x_3^6.
\end{equation}
Thus, we can assume that $g$ lies in the polynomial ring
$\R[\theta_1, \theta_2, \theta_3]$.

The first condition of Theorem~\ref{thm:cohn-elkies}
is a simple linear condition in the coefficients of the polynomial
$g$, it says
\[
g(0) \geq 1.
\]
%where we for instance have
%\[
%\vol  B^p_3 = \frac{8 \Gamma(1+1/p)^3}{\Gamma(1+3/p)}.
%\]

For the second condition we want the function $f$ whose Fourier
transform is given by \eqref{eq:fhat} to be of positive type. This is
true if and only if $g$ is globally nonnegative.  In general, 
checking that a polynomial is nonnegative everywhere is
computationally difficult, it is an $\textrm{NP}$-hard problem. We
use a standard relaxation of the global nonnegativity constraint by
imposing a sufficient condition which is easier to check: We want that
$g$ can be written as a sum of squares, which we can formulate as a
semidefinite condition. Furthermore, we can use the imposed
$\mathsf{B}_3$-invariance of $g$ to simplify this semidefinite
condition. In Section~\ref{sec:sos-pseudo-reflections} we work out the
theory of this simplification for the case of pseudo-reflection
groups. In Section~\ref{sec:sos-b3} we apply the theory
to the finite reflection group~$\mathsf{B}_3$.

Although using this sum of squares relaxation works very well in practice,
we are indeed restricting the search space of functions. Hilbert
showed in 1888 that there are polynomials already in two variables which are
globally nonnegative but which are not sum of squares. Hilbert's
proof was nonconstructive and only in 1967 Motzkin published the
first explicit example. Shortly afterwards Robinson showed that the
$\mathsf{B}_3$-invariant polynomial
\begin{equation}
\label{eq:robinson}
x_1^6 + x_2^6 + x_3^6 - (x_1^4 x_2^2 +x_1^2 x_2^4 +
x_1^4 x_3^2 + x_1^2 x_3^4 + x_2^4 x_3^2 + x_2^2 x_3^4) + 3 x_1^2 x_2^2 x_3^2
\end{equation}
is nonnegative but not a sum of squares. We refer the interested
reader to Reznick~\cite{Reznick2000a} for more on this.

For the third condition we first have to compute $f$ from
$\widehat{f}$. This is an easy linear algebra computation once we
decompose $g$ as a sum of products of radial polynomials times
harmonic polynomials. We review this decomposition in
Section~\ref{sec:fourier-transform}.

Finally, we want that $f$ be nonpositive outside
of~$\mathcal{K}^\circ - \mathcal{K}^\circ$. When $\mathcal{K} = B^p_3$ and when $p$ is an even
integer we can use another sufficient sum of squares condition:
\begin{equation}
\label{eq:sos-condition}
f(x_1,x_2,x_3) e^{\pi\|x\|^2} + (x_1^p + x_2^p  + x_3^p - 2) q_1(x_1,x_2,x_3) + q_2(x_1,x_2,x_3)  = 0
\end{equation}
where $q_1$ and~$q_2$ are $\mathsf{B}_3$-invariant polynomials which
can be written as sum of squares.  This again can be expressed as a
semidefinite condition.

So in the end we can find a good function $f$, minimizing $f(0)$, for
Theorem~\ref{thm:cohn-elkies} by solving a finite semidefinite
programming problem, once we restrict the degrees of polynomials $g$,
$q_1$, and $q_2$. We give an explicit finite-dimensional semidefinite
programming formulation in Section~\ref{sec:formulation}.

A semidefinite programming problem --- a rich generalization of linear
programming --- amounts to minimizing a linear function over an
spectrahedron, the intersection of the cone of positive semidefinite
matrices with an affine subspace. For solving semidefinite programming
problems in practice one uses interior point methods. There are many
very good software implementations of interior point methods
available. For verifying that we proved a rigorous bound we only have
to show that the solution the software gave to us is a function $f$
which satisfies the conditions of Theorem~\ref{thm:cohn-elkies}. In
Section~\ref{sec:verification} we explain this verification process in
detail.

When $p$ is not an even integer, the approach of using sum of squares
in \eqref{eq:sos-condition} breaks down. To get an upper bound we use
a sum of squares condition for the next largest even integer $p'$ and
use a fine sample of points in the intersection of the set
$B^p_3 \setminus B^{p'}_3$ with the fundamental domain
\begin{equation}
\label{eq:fundamentaldomain}
0 \leq x_1 \leq x_2 \leq x_3
\end{equation}
of the group $\mathsf{B}_3$ to make sure that the function $f$ is
nonpositive there. With this we get a function $f$ which almost
satisfies the conditions of Theorem~\ref{thm:cohn-elkies}. It turns
out, and we check this fact rigorously, that $f$ satisfies the
conditions for a slightly larger body $\alpha \mathcal{K}$ with
$\alpha$ only slightly larger than one. Then we obtain the slightly
weaker bound of $\alpha^3 f(0) \vol \mathcal{K}$.

When dealing with polytopes $\mathcal{K}$ we use a similar approach:
We impose the sum of squares condition
\begin{equation}
f(x_1,x_2,x_3) e^{\pi\|x\|^2} + (x_1^2 + x_2^2  + x_3^2 - r) q_1(x_1,x_2,x_3) + q_2(x_1,x_2,x_3)  = 0,
\end{equation}
where $r$ is the circumradius of the polytope $\mathcal{K}$.  Again we
use a fine sample of points in the intersection of the set
$rB^2_3 \setminus \mathcal{K}$ with the fundamental
domain~\eqref{eq:fundamentaldomain} of the finite reflection group
$\mathsf{B}_3$ to make sure that the function $f$ is nonpositive
there.

\subsection{Future research}

We end the introduction by showing directions and questions for
possible future research.

Our bounds give hope that the Cohn-Elkies bound might be strong enough
to prove optimality of the $C_1$-lattices for some values of
$p$ among all translative packings of superballs. Our computations
were restricted, due to numerical difficulties, to polynomials $g$ of
rather small degrees and to sum of square certificates for
nonnegativity. Does there exist a threshold $p' < \infty$ so that for all $p
\geq p'$ the Cohn-Elkies bound is tight? 

So the development of better computational techniques to compute the
Cohn-Elkies bound would be very valuable. It also would be of interest
to perform more computations. P\"utz, in his master's thesis
\cite{Puetz2016a}, computed bounds for Platonic and Archimedean solids
with icosahedral symmetries. Bounds for superball or polytope packings
in dimension~$4$ have not yet been computed.

When computing the bound for translative packings of $B^p_3$ with odd
$p$ or for translative packings of polytopes we used sampling. This
makes finding a rigorous proof more difficult. Is it possible to find a
more convenient method to prove rigorous (and better) bounds?

Is it possible to apply Minkowski's method, or a variant of the
algorithm of Betke and Henk to determine optimal lattice packings of
three-dimensional superballs?

Cohn and Zhao \cite{Cohn2014a} improve the asymptotic sphere packing
bound by Kabatiansky and Levenshtein \cite{Kabatiansky1978a} slightly
and show that the Cohn-Elkies bound is at least as strong as the
Kabatiansky-Levenshtein bound. Elkies, Odlyzko, and Rush
\cite{Elkies1991a} improve the Minkowski-Hlawka lower bound for lattice
packings of superballs. Fejes T\'oth, Fodor, and V\'\i{}gh
\cite{FejesToth2015a} find upper bounds for congruent packings of
$n$-dimensional regular cross polytopes when $n \geq 7$. How does the
Cohn-Elkies bound behave asymptotically for translative superball
packings?

With a generalization of the Cohn-Elkies bound one can also consider
packings of congruent copies of a given body, but this is
computationally even more challenging. This basic setup is explained
in Oliveira and Vallentin \cite{Oliveira2013a} where they consider
packings of congruent copies of regular pentagons in the Euclidean
plane.

\section{Sums of Hermitian squares invariant under a finite group
  generated by pseudo-reflections}
\label{sec:sos-pseudo-reflections}

Testing that a given real, multivariate polynomial is a sum of squares (SOS)
is a fundamental computational problem in polynomial optimization and
real algebraic geometry; see the recent book edited by Blekherman, Parrilo,
and Thomas \cite{Blekherman2013a}. 

Using the Gram matrix method this test can be reduced to the feasibility
problem of semidefinite optimization: A real, multivariate polynomial
$p \in \mathbb{R}[x_1, \ldots, x_n]$ of degree $2d$ is an SOS if and
only if there is a positive semidefinite matrix $Q$ of size
$\binom{n+d}{d} \times \binom{n+d}{d}$ --- a Gram matrix
representation of $p$ --- so that
\begin{equation}
\label{equality}
p(x_1, \ldots, x_n) = b(x_1, \ldots, x_n)^{\sf T} \, Q \, b(x_1, \ldots, x_n)
\end{equation}
holds, where
$b(x_1, \ldots, x_n) \in \mathbb{R}[x_1, \ldots,
x_n]^{\binom{n+d}{d}}$
is a vector which contains a basis of the space of real polynomials up
to degree $d$.

Gatermann and Parrilo~\cite{Gatermann2004a} developed a general theory
to simplify the matrices occurring in the Gram matrix method when the
polynomial at hand is invariant under the action of a finite matrix
group; see also Bachoc et al. \cite{Bachoc2012a}.

In this section we work out this simplification for polynomials
invariant under a finite group generated by pseudo-reflections. A
pseudo-reflection is a linear transformation of $\mathbb{C}^n$ where
precisely one eigenvalue is not equal to one. In particular, a
reflection $x \mapsto x - \frac{2 x \cdot v}{v \cdot v} v$ at a linear
hyperplane orthogonal to a vector $v$ is a pseudo-reflection.

In this case the computations required to apply the general theory of
Gatermann and Parrilo can be done rather concretely on the basis of
the theory developed by Shephard and Todd, Chevalley, and Serre (see
for example the book by Humphreys \cite{Humphreys1992a}, the survey by
Stanley \cite{Stanley1979a}, or the book by Sturmfels
\cite{Sturmfels1993a}).

However, we deviate from the path set out by Gatermann and Parrilo in
one important detail. Gatermann and Parrilo consider polynomials over
the field of real numbers. When working with finite groups generated
by pseudo-reflections it is more natural to work in the framework of
Hermitian symmetric polynomials since we will use the Peter-Weyl
theorem, the decomposition of the regular representation into
irreducible unitary representations.

A polynomial
$p \in \C[z_1, \ldots, z_n, \overline{w}_1, \ldots, \overline{w}_n] =
\C[z,\overline{w}]$ is called a \defi{Hermitian symmetric polynomial} if one
of the following three equivalent conditions holds (see D'Angelo \cite{DAngelo2011a}):
\begin{enumerate}
\item[i)] Equality $p(z,\overline{w}) = \overline{p(w,\overline{z})}$
  holds for all $z, w \in \C^n$.
\item[ii)] The function $z \mapsto p(z,\overline{z})$, with $z \in
  \mathbb{C}^n$, is real-valued.
\item[iii)] There is a Hermitian matrix $Q = (q_{\alpha \beta})$ so
  that one can represent $p$ as
  $p(z,\overline{w}) = \sum_{\alpha,\beta} q_{\alpha \beta} z^\alpha
  \overline{w}^\beta$.
\end{enumerate}

A Hermitian symmetric polynomial $p \in \C[z,\overline{w}]$ is a
\defi{sum of Hermitian squares} if there are polynomials
$q_1, \ldots, q_r \in \C[z]$ so that
\[
p(z,\overline{w}) = \sum_{i=1}^r q_i(z) \overline{q_i(w)}
\]
holds. In particular, a sum of Hermitian squares determines a
real-valued nonnegative function by
$z \mapsto p(z,\overline{z}) = \sum_{i=1}^r |q_i(z)|^2$. In fact,
D'Angelo gave in \cite[Definition IV.5.1]{DAngelo2002a} eight
positivity conditions for a Hermitian symmetric polynomial. He noted
that being a sum of Hermitian squares is the strongest among them and
that this condition is also easy to verify by the Gram matrix method
after an obvious adaptation of \eqref{equality}:
\[
p(z, \overline{w}) = b(z)^{\sf T} \, Q \, \overline{b(w)}, 
\]
where $Q$ is a Hermitian positive semidefinite matrix of size
$\binom{n+d}{d} \times \binom{n+d}{d}$ and where
$b(z) \in \mathbb{C}[z]^{\binom{n+d}{d}}$
is a vector which contains a basis of the space of complex polynomials
up to degree $d$.

Now let us review the relevant theory of pseudo-reflection groups: Let
$G \subseteq \mathrm{GL}_n(\mathbb{C})$ be a finite group generated by
pseudo-reflections. It is acting on the polynomial ring
$\mathbb{C}[x_1, \ldots, x_n] = \mathbb{C}[x]$ by
\[
(gp)(x) = p(g^{-1}x) \quad \text{for $g \in G$ and $p \in \C[x]$.}
\]
The invariant ring is defined by
\[
\mathbb{C}[x]^G = \{p \in \mathbb{C}[x] : gp = p \text{ for all } g
\in G\}.
\]
The invariant ring is generated by $n$ homogeneous polynomials $\theta_1, \ldots,
\theta_n$ which are algebraically independent. Thus,
\[
\mathbb{C}[x]^G =\mathbb{C}[\theta_1, \ldots, \theta_n]
\]
is a free algebra. Homogeneous, algebraically independent generators of
the invariant ring are called basic invariants. They are not uniquely
determined by the group, but their degrees $d_1, \ldots, d_n$ are. 

The group action respects the grading of the polynomial ring. To
determine the dimensions of the invariant subspaces of homogeneous
polynomials
\[
\Hom^G_k = \mathbb{C}[x]^G \, \cap \, \Hom_k
\]
with
\[  
\Hom_k= \left\{p \in \C[x] : p(\alpha x) = \alpha^k p(x) \, \text{for
    all } \alpha \in \C, \deg p = k\right\},
\]
one can use Molien's series
\[
\sum_{k = 0}^\infty \dim \Hom^G_k t^k = \left(\prod_{i=1}^n (1-t^{d_i})\right)^{-1}.
\]

The coinvariant algebra is
\[
\mathbb{C}[x]_G = \mathbb{C}[x] / I, 
\]
where $I = (\theta_1, \ldots, \theta_n)$ is the ideal generated by
basic invariants. The coinvariant algebra is a graded algebra of
finite dimension $|G|$.  The dimensions of the homogeneous subspaces
of $\mathbb{C}[x]_G$ are given by the Poincar\'e series
\[
(1-t)^{-n} \prod_{i=1}^n (1-t^{d_i}).
\]

In particular,
\[
\mathbb{C}[x] = \mathbb{C}[x]^G \otimes \mathbb{C}[x]_G
\]
holds. The action of $G$ on the coinvariant algebra $\mathbb{C}[x]_G$ is
equivalent to the regular representation of $G$. Let $\widehat{G}$ be
the set of irreducible unitary representations of $G$ up to
equivalence. Then one can apply the Peter-Weyl theorem, see for example \cite[Chapter 15]{Terras1999a}: There are
homogeneous polynomials
\[
\varphi^{\pi}_{ij}, \; \text{ with } \; \pi \in \widehat{G}, \; 1 \leq
i,j \leq d_{\pi},
\]
where $d_{\pi}$ is the degree of $\pi$, which form a basis of the
coinvariant algebra such that the transformation law
\begin{equation}
\label{eq:transformation}
g\varphi^{\pi}_{ij} = (\pi(g)_j)^{\sf T}
\begin{pmatrix}
\varphi^{\pi}_{i1}\\
\vdots\\
\varphi^{\pi}_{id_{\pi}}
\end{pmatrix},
\;
i = 1, \ldots, d_\pi,
\end{equation}
holds for all $g \in G$. Here, $\pi(g)_j$ denotes the $j$-th column of
the unitary matrix $\pi(g) \in \mathrm{U}(d_\pi)$.

We extend the action of $G$ from $\mathbb{C}[x]$ to
$\mathbb{C}[z, \overline{w}]$ by
\[
(gp)(z,\overline{w}) = p(g^{-1}z,\overline{g^{-1}w}) \quad \text{for $g \in G$ and $p \in \C[z,\overline{w}]$.}
\]
We define the ring of $G$-invariant Hermitian symmetric
polynomials by
\[
\mathbb{C}[z,\overline{w}]^G = \{p \in \mathbb{C}[z,\overline{w}] : gp = p \text{ for all } g
\in G\}.
\]

Now we set up all necessary notation for formulating the theorem which
gives an explicit parametrization of the convex cone of $G$-invariant
Hermitian symmetric polynomials which are Hermitian sum of squares.
The following theorem can be derived from the real version of
\cite[Theorem 6.2]{Gatermann2004a}. So we omit the proof.

\begin{theorem}
\label{thm:invariant-sos}
 Let $G \subseteq \mathrm{GL}_n(\mathbb{C})$ be a
  finite group generated by pseudo-reflections. The convex cone of
  $G$-invariant Hermitian symmetric polynomials which can be written as sums of Hermitian squares
  equals
\[
\begin{split}
\Big\{
p \in \mathbb{C}[z,\overline{w}]^G : \; &  p(z,\overline{w}) = \sum_{\pi
  \in \widehat{G}} \langle P^{\pi}(z,\overline{w}),
Q^{\pi}(z,\overline{w}) \rangle,\\
& \text{$P^{\pi}(z,\overline{w})$ is a Hermitian SOS matrix
polynomial in $\theta_i$} \Big\}.
\end{split}
\]
Here $\langle A, B \rangle = \mathrm{Tr}(B^*A)$ denotes the trace
inner product, the matrix $P^{\pi}(z,\overline{w})$ is a Hermitian SOS
matrix polynomial in the variables $\theta_1, \ldots, \theta_n$, i.e.\
there is a matrix $L^{\pi}(z)$ with entries in
$\mathbb{C}[z]^G = \mathbb{C}[\theta_1, \ldots, \theta_n]$ such that
\[
P^\pi(z,\overline{w}) = L^{\pi}(z) \overline{L^{\pi}(w)}^{\sf T} 
\]
holds, and
$Q^{\pi}(z,\overline{w}) \in (\mathbb{C}[z,\overline{w}]^G)^{d_\pi
  \times d_\pi}$ is defined componentwise by
\[
[Q^\pi]_{kl}(z,\overline{w}) = \sum_{i=1}^{d_\pi} \varphi^{\pi}_{ki}(z) \overline{\varphi^{\pi}_{li}(w)}.
\]
\end{theorem}

The computational value of this approach is that one only has to
determine basic invariants $\theta_1, \ldots, \theta_n$ and a suitable
basis $\varphi^{\pi}_{ij}$ of the coinvariant algebra which
satisfies~\eqref{eq:transformation}. These computations are
\textit{independent} of the degree of the polynomial~$p$.

\medskip

It turns out that for the octahedral group $\mathsf{B}_3$ we consider
for our application all irreducible unitary representation are
orthogonal representations. In this case the previous theorem can be
translated into the following version for the field of real numbers.

\begin{theorem}
\label{thm:real-invariant-sos}
Let $G \subseteq \mathrm{GL}_n(\mathbb{R})$ be a finite group
generated by pseudo-reflections so that all unitary irreducible
representation $\pi \in \widehat{G}$ of $G$ are orthogonal. The convex
cone of $G$-invariant real polynomials which can be written as sums of
squares equals
\[
\begin{split}
\Big\{
p \in \mathbb{R}[x]^G : \; &  p(x) = \sum_{\pi
  \in \widehat{G}} \langle P^{\pi}(x),
Q^{\pi}(x) \rangle,\\
& \text{$P^{\pi}(x)$  is an SOS matrix
polynomial in $\theta_i$} \Big\}.
\end{split}
\]
Here the matrix $P^{\pi}(x)$ is an SOS matrix
polynomial in the variables $\theta_1, \ldots, \theta_n$, i.e.\ there
is a matrix $L^{\pi}(x)$ with entries in $\mathbb{R}[x]^G =
\mathbb{R}[\theta_1, \ldots, \theta_n]$ such that
\[
P^\pi(x) = L^{\pi}(x) L^{\pi}(x)^{\sf T} 
\]
holds, and
$Q^{\pi}(x) \in (\mathbb{R}[x]^G)^{d_\pi
  \times d_\pi}$ is defined componentwise by
\[
[Q^\pi]_{kl}(x) = \sum_{i=1}^{d_\pi} \varphi^{\pi}_{ki}(x) \varphi^{\pi}_{li}(x).
\]
\end{theorem}

\section{Real sums of squares polynomials invariant under the
  octahedral group}
\label{sec:sos-b3}

In this section we specialize Theorem~\ref{thm:real-invariant-sos} to
the symmetry group of the three-dimensional real octahedron, the
octahedral group, which is the finite reflection group $\mathsf{B}_3$
generated by the matrices~\eqref{eq:b3-generators}. Since in the
literature only very few cases of Theorem~\ref{thm:invariant-sos} or
Theorem~\ref{thm:real-invariant-sos} are worked out explicitly, we
give substantial amount of detail here.

We use the basic invariants $\theta_1, \theta_2, \theta_3$ as given in
\eqref{eq:basic-invariants}. Let $\chi_{\pi}$ be the character of the
irreducible representations $\pi \in \widehat{\mathsf{B}}_3$. There
are ten inequivalent irreducible unitary representations and the
character table, which one computes with a computer algebra system or
which one also can find in many text books on mathematical chemistry,
is given in Table~\ref{table:character-table}.

% nice diagrams: http://kociemba.org/math/symmetric.htm
% http://gernot-katzers-spice-pages.com/character_tables/Oh.html

\begin{table}[htb]

\renewcommand{\arraystretch}{1.2}

\begin{center}
\begin{tabular}{c|*{10}{r}}
\hline
                & $E$ & $i$  & $3C_2$ & $3\sigma_h$ & $6C_2'$ & $6\sigma_d$ & $8C_3$ & $6C_4$ & $6 S_4$ & $8S_6$ \\
\hline
$A_{1g}$ & $1$ & $ 1$ & $1$       & $1$                &$1$         & $1$                & $1$       & $1$        & $1$       & $1$ \\ 
$A_{1u}$ & $1$ & $-1$ & $1$       & $-1$               &$1$         & $-1$              & $1$        & $1$       & $-1$       & $-1$ \\
$A_{2g}$ & $1$ & $ 1$ & $1$       & $1$                & $-1$        & $-1$              & $1$       & $-1$       & $-1$       & $1$ \\
$A_{2u}$ & $1$ & $-1$ & $1$       & $-1$              &$-1$         & $1$               & $1$        & $-1$       & $1$      & $-1$ \\
$E_g$     & $2$ & $ 2$ & $2$        & $2$               &$0$         & $0$                & $-1$        & $0$      & $0$      & $-1$ \\
$E_u$     & $2$ & $-2$ & $2$        & $-2$             & $0$         & $0$               & $-1$        & $0$      & $0$       & $1$ \\
$T_{1g}$ & $3$ & $ 3$ & $-1$      & $-1$              & $-1$        & $-1$             & $0$        & $1$       & $1$       & $0$ \\
$T_{1u}$ & $3$ & $-3$ & $-1$      & $1$               & $-1$        & $1$              & $0$        & $1$       & $-1$       & $0$ \\
$T_{2g}$ & $3$ & $ 3$ & $-1$      & $-1$             & $1$        & $1$                & $0$        & $-1$       & $-1$       & $0$ \\
$T_{2u}$ & $3$ & $-3$ & $-1$      & $1$              & $1$        & $-1$               & $0$        & $-1$      & $1$        & $0$ \\
\hline
\end{tabular}
\\[0.3cm]
\end{center}
\caption{Character table of $\mathsf{B}_3$. Rows are indexed by characters, and columns are indexed by conjugacy classes.}
\label{table:character-table}
\end{table}

In the character table we use Mulliken symbols for concreteness. This
scheme was suggested by Robert S.~Mulliken, Nobel laureate in
Chemistry in 1966. The symmetry group of the regular three-dimensional
octahedron coincides with the one of the regular cube. In the
following we describe the conjugacy classes of $\mathsf{B}_3$
geometrically by looking at the symmetries of the cube: $E$ is the
identity of the group ($E$ from German \textit{Einheit}), $i$ is the
inverse operation $i(x) = -x$, $3C_2$ are the three clockwise
rotations by $180^\circ$ through the axis of the facet centers,
$3\sigma_h$ are the three reflections through planes which are
parallel to pairs of facets ($\sigma$ from \textit{Spiegelung}),
$6C_2'$ are the six clockwise rotations by $180^\circ$ through the
axis of the edge centers, $6\sigma_d$ are the six reflections through
the planes given by the diagonals of the facets, $8C_3$ are the eight
clockwise rotations by $120^\circ$ through the diagonals of the cube,
$6C_4$ are the six clockwise rotations by $90^\circ$ through the axis
of the facet centers, $6 S_4$ are the six rotation-reflections by
$90^\circ$ through the axis of the facet centers, and $8S_6$ are the
eight rotation-reflections by $60^\circ$ through the diagonals of the
cube. One-dimensional characters are given by the letter $A$,
two-dimensional characters are specified by the letter $E$, and the
three-dimensional ones by $T$. The subscript $g$ (\textit{gerade}) or
$u$ (\textit{ungerade}) is used to distinguish between $\chi(i) = 1$
and $\chi(i) = -1$.

The Molien series of $\C[x_1,x_2,x_3]^{\mathsf{B}_3}$ is
\[
\frac{1}{(1-t^2)(1-t^4)(1-t^6)} = 1 + t^2 + 2t^4 + 3t^6 + 4t^8 +
5t^{10} + 7t^{12} + 8t^{14} + 10t^{16} + 12t^{18} + \cdots
\]
The coinvariant algebra $\C[x]_G = \C[x]/I$ with
$I = (\theta_1, \theta_2, \theta_3)$ decomposes into
\[
V = V_0 \oplus \cdots \oplus V_9
\]
according to the grading by degree, where the dimensions of the spaces $V_k$,
with $k = 0, \ldots, 9$ can be read off by the Poincar\'e series
\[
1 + 3t + 5t^2 + 7t^3 + 8t^4 + 8t^5 + 7t^6 + 5t^7 + 3t^8 + t^9.
\]
The group action respects the grading. It turns out that all
irreducible unitary representations occur multiplicity-free in the
$V_k$'s and that all of them are orthogonal representations. Serre
\cite[Chapter 2.6, Theorem 8]{Serre1977a} gives a formula which can be
used to decompose a finite-dimensional representation into its
isotypic components. Consider the representation
\[
\rho_k \colon \mathsf{B}_3 \to \GL(\Hom_k),\quad \rho_k(g)(p) \mapsto gp,
\]
and consider a unitary irreducible representation
$\pi \in \widehat{\mathsf{B}}_3$. Then the image of the linear map
\[
p^\pi_k \colon \Hom_k \to \Hom_k, \quad p^\pi_k = \frac{d_{\pi}}{|\mathsf{B}_3|} \sum_{g \in \mathsf{B}_3} \chi_\pi(g^{-1}) \rho_k(g)
\]
gives the subspace $V^\pi_k$ of $\Hom_k$ which is the isotypic
component of $\Hom_k$ having type $\pi$.

We choose the smallest degree $k_\pi$ so that there is a nontrivial
isotypic component of $\Hom_{k_\pi}$ having type $\pi$, and this
choice of~$k_\pi$ implies that this isotypic component $V^\pi_{k_\pi}$ is actually an
irreducible subspace. Then we equip this irreducible subspace with a
$\mathsf{B}_3$-invariant inner product and compute an orthonormal
basis by Gram-Schmidt orthonormalization. This orthonormal basis gives
polynomials $\varphi^{\pi}_{1j}$, with $j = 1, \ldots, d_{\pi}$, which
we need for applying Theorem~\ref{thm:real-invariant-sos}. The results
are displayed in Table~\ref{table:onb}.

\begin{table}[htb]
\begin{center}
\renewcommand{\arraystretch}{1.2}
\begin{tabular}{c|l}
\hline
$A_{1g}$ & $1$   \\
$A_{1u}$ &  $x_1x_2x_3$ \\
$A_{2g}$ &  $x_1^4x_2^2 - x_1^4x_3^2 - x_1^2x_2^4 + x_1^2x_3^4 + x_2^4x_3^2 - x_2^2x_3^4$\\
$A_{2u}$ & $x_1^5x_2^3x_3 - x_1^5x_2x_3^3 - x_1^3x_2^5x_3 + x_1^3x_2x_3^5 + x_1x_2^5x_3^3 - x_1x_2^3x_3^5$\\
$E_g$     & $x_1^3x_2x_3 - x_1x_2x_3^3$\\
              & $\frac{\sqrt{3}}{3} x_1^3x_2x_3 - \frac{2\sqrt{3}}{3} x_1x_2^3x_3
                + \frac{\sqrt{3}}{3} x_1x_2x_3^3$\\
$E_u$     &  $x_1^2 - x_3^2$\\
               & $\frac{\sqrt{3}}{3}x_1^2 - \frac{2 \sqrt{3}}{3} x_2^2 + \frac{\sqrt{3}}{3}x_3^2$\\
$T_{1g}$ & $x_1^3x_2 - x_1x_2^3$\\
             & $x_1^3x_3 - x_1x_3^3$\\
             & $x_2^3x_3 - x_2x_3^3$ \\
$T_{1u}$ & $x_1^2x_2 - x_2x_3^2$\\
             & $x_1^2x_3 - x_2^2x_3$\\
             & $x_1x_2^2 - x_1x_3^2$ \\
$T_{2g}$ & $x_1x_2$\\
              & $x_1x_3$\\
              & $x_2x_3$\\
$T_{2u}$ & $x_1$   \\
               & $x_2$  \\
               & $x_3$  \\
\hline
\end{tabular}
\\[0.3cm]
\end{center}
\caption{Orthonormal basis $\varphi^\pi_{1j}$, with $j = 1, \ldots,
  d_\pi$, of subspaces $V^\pi_{k_\pi}$.}
\label{table:onb}
\end{table}

The next task is to find the other polynomials $\varphi^\pi_{ij}$,
with $i = 2, \ldots, d_{\pi}$, which transform according to
\eqref{eq:transformation}. We use the algorithm of Serre \cite[Chapter
2.7, Proposition 8]{Serre1977a} for this.

Define
\[
p^{\pi}_{k,ij} \colon \Hom_k \to \Hom_k \quad \text{by} \quad p^{\pi}_{k,ij} = \frac{d_{\pi}}{|\mathsf{B}_3|} \sum_{g \in \mathsf{B}_3} \pi_{ji}(g^{-1}) \rho_k(g),
\]
where $\pi(g) \in \mathrm{U}(d_\pi)$ is the unitary matrix which we
get by considering the matrix representation $\rho_{k_{\pi}}(g)$
restricted to the irreducible subspace of $V^\pi_{k_\pi}$ of $\Hom_{k_{\pi}}$ having type
$\pi$ and expressed in terms of the orthonormal basis
$\varphi^\pi_{1j}$, with $j = 1, \ldots, d_\pi$, we just
computed. Denote the image $p^{\pi}_{k,ii}(V^{\pi}_k)$ by
$V^{\pi}_{k,i}$, where $i = 1, \ldots, d_\pi$. Then we have the
decomposition
\[
V^{\pi}_k = V^\pi_{k,1} \oplus \cdots \oplus V^\pi_{k,d_\pi}.
\]
Consider a nonzero vector $\varphi^{\pi}_{k1} \in
V^\pi_{k,1}$. Define $\varphi^{\pi}_{ki} = p^\pi_{k,i1}(\varphi^\pi_{k1})$. Then,
\[
\rho_k(g)(\varphi^{\pi}_{ki}) = \sum_{j=1}^{d_\pi} \pi_{ji}(g) \varphi^{\pi}_{kj}
\]
holds for all $g \in G$, as we wanted. With this information we can
construct the matrices $Q^{\pi}$. We give them in
Table~\ref{table:qpi}.

\begin{table}
\begin{center}
\renewcommand{\arraystretch}{1.2}
\begin{tabular}{c|l}
\hline
$A_{1g}$ & $1$\\
$A_{1u}$ & $\theta_1^3 - 3\theta_1\theta_2 + 2\theta_3$\\
$A_{2g}$ & $-\theta_1^6 + 9\theta_1^4\theta_2 - 8\theta_1^3\theta_3 - 21\theta_1^2\theta_2^2 + 36\theta_1\theta_2\theta_3 + 3\theta_2^3 - 18\theta_3^2$\\
$A_{2u}$ & $-\theta_1^9 + 12\theta_1^7\theta_2 - 10\theta_1^6\theta_3 -
48\theta_1^5\theta_2^2 + 78\theta_1^4\theta_2\theta_3 +
66\theta_1^3\theta_2^3 - 34\theta_1^3\theta_3^2 -
150\theta_1^2\theta_2^2\theta_3$ \\
& \qquad $- 9\theta_1\theta_2^4 + 126\theta_1\theta_2\theta_3^2 +
  6\theta_2^3\theta_3 - 36\theta_3^3$\\
$E_g$ & $-2\theta_1^5 + 12\theta_1^3\theta_2 - 4\theta_1^2\theta_3 -
18\theta_1\theta_2^2 + 12\theta_2\theta_3$\\
&  $-2\theta_1^4\theta_2 +
6\theta_1^3\theta_3 + 6\theta_1^2\theta_2^2 -
22\theta_1\theta_2\theta_3 + 12\theta_3^2$\\
&  $\theta_1^7 -
9\theta_1^5\theta_2 + 10\theta_1^4\theta_3 + 19\theta_1^3\theta_2^2 -
36\theta_1^2\theta_2\theta_3 - 3\theta_1\theta_2^3 +
16\theta_1\theta_3^2 + 2\theta_2^2\theta_3$\\
$E_u$ & $-2\theta_1^2 + 6\theta_2$\\
& $-2\theta_1\theta_2 + 6\theta_3$\\
& $\theta_1^4 - 6\theta_1^2\theta_2 + 8\theta_1\theta_3 +
  \theta_2^2$\\
$T_{1g}$ & $12\theta_1\theta_3 - 12\theta_2^2$\\
& $2\theta_1^5 - 12\theta_1^3\theta_2 + 16\theta_1^2\theta_3 +
6\theta_1\theta_2^2 - 12\theta_2\theta_3$\\
&
$2\theta_1^6 - 12\theta_1^4\theta_2 + 10\theta_1^3\theta_3 +
12\theta_1^2\theta_2^2 - 6\theta_1\theta_2\theta_3 - 6\theta_2^3$\\
&
$2\theta_1^6 - 10\theta_1^4\theta_2 + 10\theta_1^3\theta_3 +
10\theta_1\theta_2\theta_3 - 12\theta_3^2$\\
&
$\theta_1^7 - 3\theta_1^5\theta_2 + 2\theta_1^4\theta_3 -
  9\theta_1^3\theta_2^2 + 24\theta_1^2\theta_2\theta_3 + 
    3\theta_1\theta_2^3 - 12\theta_1\theta_3^2 -
  6\theta_2^2\theta_3$\\
&
$4\theta_1^6\theta_2 -
    3\theta_1^5\theta_3 - 21\theta_1^4\theta_2^2 + 32\theta_1^3\theta_2\theta_3 + 12\theta_1^2\theta_2^3 - 12\theta_1^2\theta_3^2 -
9\theta_1\theta_2^2\theta_3 - 3\theta_2^4$\\
$T_{1u}$ & $-12\theta_1^3 + 48\theta_1\theta_2 - 36\theta_3$\\
& $-6\theta_1^4 + 24\theta_1^2\theta_2 - 12\theta_1\theta_3 - 6\theta_2^2$\\
& $-6\theta_1^3\theta_2 + 6\theta_1^2\theta_3 + 18\theta_1\theta_2^2 - 18\theta_2\theta_3$\\
& $-2\theta_1^5 + 6\theta_1^3\theta_2 + 2\theta_1^2\theta_3 - 6\theta_2\theta_3$\\
& $\theta_1^6 - 9\theta_1^4\theta_2 + 8\theta_1^3\theta_3 + 15\theta_1^2\theta_2^2 - 12\theta_1\theta_2\theta_3 - 3\theta_2^3$\\
& $\theta_1^7 - 6\theta_1^5\theta_2 + 5\theta_1^4\theta_3 + 3\theta_1^3\theta_2^2 + 6\theta_1\theta_2^3 - 9\theta_2^2\theta_3$\\
$T_{2g}$ & $3\theta_1^2 - 3\theta_2$\\
& $6\theta_1\theta_2 - 6\theta_3$\\
& $-\theta_1^4 + 6\theta_1^2\theta_2 - 2\theta_1\theta_3 - 3\theta_2^2$\\
& $-2\theta_1^4 + 12\theta_1^2\theta_2 - 10\theta_1\theta_3$\\
& $-\theta_1^5 + 4\theta_1^3\theta_2 - 2\theta_1^2\theta_3 + 3\theta_1\theta_2^2 - 4\theta_2\theta_3$\\
& $-2\theta_1^4\theta_2 + \theta_1^3\theta_3 + 9\theta_1^2\theta_2^2 - 7\theta_1\theta_2\theta_3 - 3\theta_2^3 + 2\theta_3^2$\\
$T_{2u}$ & $6\theta_1$\\
& $6\theta_2$\\
& $6\theta_3$\\
& $6\theta_3$\\
& $\theta_1^4 - 6\theta_1^2\theta_2 + 8\theta_1\theta_3 + 3\theta_2^2$\\
& $\theta_1^5 - 5\theta_1^3\theta_2 + 5\theta_1^2\theta_3 + 5\theta_2\theta_3$\\
\hline
\end{tabular}
\\[0.3cm]
\end{center}
\caption{Matrices $Q^{\pi}$ for the group $\mathsf{B}_3$ given in
  upper triangular row-major order (in the consecutive order of row entries of
  upper triangular matrices).}
\label{table:qpi}
\end{table}

\section{Computing the Fourier transform}
\label{sec:fourier-transform}

As explained in the introduction we define the function $f$ which we
want to use in Theorem~\ref{thm:cohn-elkies} through its Fourier
transform $\widehat{f}(u) = g(u) e^{-\pi \|u\|^2}$ where $g$ is a
polynomial. In order to verify the third condition of the theorem, we
have to compute $f$ from $\widehat{f}$. In other words, we have to
compute the Fourier transform of $u \mapsto \widehat{f}(-u)$. In this
section we explain how to do this. We first consider the general case,
when $g$ is an arbitrary complex polynomial in $n$ variables. Then we
show how some of the computations can be simplified when we assume
that $g$ is $\mathsf{B}_3$-invariant. In the end, since the Fourier
transform is linear, we have to solve a certain system of linear
equations.  A similar calculation was done by Dunkl~\cite{Dunkl1981a}.
He even gives explicit algebraic solutions.

Consider the following decomposition of complex polynomials in $n$
variables of degree at most $d$:
\begin{equation}
\label{eq:decomposition}
\C[x]_{\leq d} = \bigoplus_{j=0}^d \Hom_j = \bigoplus_{j=0}^d \bigoplus_{\substack{r,k\\ 2r+k = d}} \|x\|^{2r} \Harm_k,
\end{equation}
where
\[
\Harm_k = \left\{h \in \Hom_k : \Delta h = \left(\frac{\partial^2}{\partial x_1^2} + \cdots + \frac{\partial^2}{\partial x_n^2}\right) h = 0\right\}
\]
is the space of (homogeneous) harmonic polynomials of degree $k$. In
other words, harmonic polynomials of degree $k$ are the kernel of the
Laplace operator
\[
\Harm_k = \ker \Delta, \quad \Delta : \Hom_k \to \Hom_{k-2},
\]
where
\[
\dim \Harm_k = \dim \Hom_k - \dim \Hom_{k-2} = \binom{n+k-1}{k} - \binom{n+k-3}{k-2}.
\]
The existence of decomposition~\eqref{eq:decomposition} is classical;
one can find a proof, for example, in the book by Stein and Weiss
\cite[Theorem IV.2.10]{Stein1971a}. Decomposition~\eqref{eq:decomposition} together with
the following proposition shows how to compute $f$ from $g$ by solving
a system of linear equations. The proposition in particular
shows that the function $x \mapsto h_k(x) e^{-\pi \|x\|^2}$ with
  $h_k \in \Harm_k$ is an eigenfunction of the Fourier transform with
  eigenvalue $i^{-k}$.

\begin{proposition}
Let
\[
f(x) = h_k(x)  \|x\|^{2r}  e^{-\pi \|x\|^2}
\]
be a Schwartz function with $h_k \in \Harm_k$. The Fourier transform
of $f$ is
\[
\widehat{f}(u) = (i^{-k} h_k(u)) \cdot \pi^{-r} r! L_r^{n/2+k-1}(\pi\|u\|^2)  e^{-\pi \|u\|^2},
\]
where $L_r^{n/2+k-1}$ is the Laguerre polynomial of degree~$r$ with
parameter~$n/2 + k - 1$.
\end{proposition}

In general, Laguerre polynomials $L^\alpha_r$ with parameter $\alpha$
are orthogonal polynomials for the inner product
$\int_0^\infty f(x) g(x) x^\alpha e^{-x} \, dx$, see the book by
Andrews, Askey, and Roy~\cite{Andrews1999a} for more details.

\begin{proof}
  Using Stein and Weiss \cite[Theorem IV.3.10]{Stein1971a} one sees
  that the Fourier transform of $f$ is
\[
\widehat{f}(u) = (i^{-k} h_k(u)) \cdot 2\pi \|u\|^{-(n/2-1+k)} \int_0^{\infty} s^{n/2+2r+k} J_{n/2-1+k}(2\pi s \|u\|)  e^{-\pi s^2} ds,
\]
where $J_{\alpha}$ is the Bessel function of the first kind of order
$\alpha$. By Andrews, Askey, and Roy \cite[Corollary 4.11.8]{Andrews1999a} the
integral above equals
\[
\frac{\Gamma(n/2+r+k) (\sqrt{\pi} \|u\|)^{n/2-1+k} e^{-\pi \|u\|^2}}{2 (\sqrt{\pi})^{n/2+2r+k+1} \Gamma(n/2+k)} {}_1 F_1\biggl({-r\atop n/2+k}; \pi \|u\|^2\biggr),
\]
where ${}_1 F_1$ denotes the hypergeometric series. Hence,
\[
\widehat{f}(u) = (i^{-k} h_k(u)) \cdot \pi^{-r} \frac{\Gamma(n/2+r+k)}{\Gamma(n/2+k)} {}_1 F_1\biggl({-r\atop n/2+k}; \pi \|u\|^2\biggr) e^{-\pi \|u\|^2}.
\]
The hypergeometric series becomes a Laguerre polynomial
(\cite[(6.2.2)]{Andrews1999a})
\[
{}_1 F_1\biggl({-r\atop n/2+k}; \pi \|u\|^2\biggr) 
= \frac{r!}{(n/2+k)_r} L_r^{n/2+k-1}(\pi\|u\|^2).
\]
Combining the last two equations gives the desired result.
\end{proof}

If one assumes that polynomial $g$ is $\mathsf{B}_3$-invariant one
can save quite some computations. Instead of working with 
decomposition~\eqref{eq:decomposition} we can work with a
$\mathsf{B}_3$-invariant decomposition because the Laplacian $\Delta$
commutes with the action of the orthogonal group:
\[
\C[x]^{\mathsf{B}_3}_{\leq d} = \bigoplus_{j=0}^d \Hom^{\mathsf{B}_3}_j = \bigoplus_{j=0}^d \bigoplus_{\substack{r,k\\ 2r+k = d}} \theta_1^r \Harm^{\mathsf{B}_3}_k.
\]
To see the computational advantage, let us compare the dimensions of
the harmonic subspaces.  We generally have $\dim \Harm_k = 2k+1$ when
$n = 3$, but the Molien series counting the dimensions of the invariant
harmonic subspaces (see Goethals and Seidel \cite{Goethals1981a}) is 
\[
\sum_{t=0}^\infty \dim \Harm^{\mathsf{B}_3}_k t^k =
\frac{1}{(1-t^4)(1-t^6)} = 1 + t^4 + t^6 + t^8 + t^{10} + 2 t^{12} +
t^{14} + 2t^{16} + 2t^{18} + \cdots .
\]

\section{Semidefinite programming formulation}
\label{sec:formulation}

We now present in detail the semidefinite program we use to find good
functions~$f$ satisfying the conditions of
Theorem~\ref{thm:cohn-elkies} when~$\Kcal \subseteq \R^3$ is such that
the Minkowski difference~$\Kcal - \Kcal$ is invariant under the action
of~$\mathsf{B}_3$. This is the case, e.g., when~$\Kcal$ is a
three-dimensional superball or a Platonic or an Archimedean solid with
tetrahedral symmetry.

\subsection{Representation of the function $f$ via its Fourier transform}

Recall that we specify the function~$f\colon \R^3 \to \R$ via its
Fourier transform. Given a real
polynomial~$g \in \R[x] = \R[x_1, x_2, x_3]$, we define
\begin{equation}
\label{eq:fhat-definition}
\widehat{f}(u) =g(u) e^{-\pi\|u\|^2}.
\end{equation}
We deal exclusively with $\bt$-invariant functions, so we take the
polynomial~$g$ above $\bt$-invariant. Functions invariant under $\bt$
are even, and so are their Fourier transforms. From this it follows
that there is no loss of generality in considering real-valued Fourier
transforms and so there is also no loss of generality in requiring
that~$g$ be a real polynomial. This simplifies the use of semidefinite
programming considerably since we only have to optimize over the cone
of real positive semidefinite matrices and not over the larger cone of
Hermitian positive semidefinite matrices.

Function~$f$ is of positive type if and only if~$g$ is a nonnegative
polynomial. Since it is computationally difficult to work with
nonnegative polynomials, we require instead that~$g$ be a sum of
squares, thus restricting (see for example the Robinson
polynomial~\eqref{eq:robinson}) the set of functions~$f$ that we work
with.

Theorem~\ref{thm:real-invariant-sos} provides a parametrization of the
cone of SOS polynomials invariant under~$\bt$, like~$g$. In the
theorem, each matrix~$P^\pi$ is an SOS matrix polynomial, that is,
there is a matrix~$L^\pi$ whose entries are invariant polynomials such
that~$P^\pi = L^\pi (L^\pi)^\tp$. To find an SOS polynomial~$g$, we
may then fix the maximum degree a polynomial in~$P^\pi$ can have, and
use the fact (cf.~Gatermann and Parrilo~\cite[Definition
2.2]{Gatermann2004a}) that~$S \in \R[x]^{n \times n}$ is an SOS matrix
if and only if the polynomial~$y^\tp S y \in \mathbb{R}[x,y]$ is a sum
of squares, where~$y = (y_1, \ldots, y_n)$ are new variables. This,
together with~\eqref{equality}, suggests a way to represent~$g$ with
one positive semidefinite matrix for each of the irreducible
representations of~$\bt$.

In our formulation we use a derived parametrization that produces
numerically stabler problems providing bounds that can be rigorously
shown to be correct. Our approach is as follows.  For each irreducible
unitary representation~$\pi \in \btd$,
let~$\Phi^\pi = \bigl(\varphi^\pi_{ij}\bigr)_{i,j=1}^{d_\pi}$ where
the~$\varphi^\pi_{ij}$'s were defined in
Section~\ref{sec:sos-pseudo-reflections} and Section~\ref{sec:sos-b3}.
Then~$Q^\pi = \Phi^\pi (\Phi^\pi)^\tp$. Each row of~$\Phi^\pi$
contains homogeneous invariant polynomials all of the same degree; we
say the \defi{degree} of a row is the degree of the polynomials in it.

Let~$\Acal$ be some basis of~$\R[x]^{\bt}$ consisting of homogeneous
polynomials. For an integer~$t \geq 0$ and each~$\pi \in \btd$,
let~$\Ical^t_\pi$ be the set of pairs~$(a, r)$, where~$a \in \Acal$
and~$1 \leq r \leq d_\pi$ indexes a row of~$\Phi^\pi$ such that the
degree of~$a$ plus the degree of the row~$r$ is at most~$t$. For
each~$\pi \in \btd$ we may then consider the matrix~$V^{\pi, t}$ with
rows and columns indexed by~$\Ical^t_\pi$ such that
\[
V^{\pi, t}_{(a, r), (b, s)} = a b Q^\pi_{rs}.
\]
Notice that entry~$((a, r), (b, s))$ of~$V^{\pi, t}$ has degree
equal to~$\deg a + \deg b + \deg Q^\pi_{rs} \leq 2t$.

In our formulation we will fix an odd\footnote{The reason why we pick
  odd~$d$ is so that the resulting problem admits a strictly feasible
  solution. This will be better explained in
  Section~\ref{sec:solving}.} positive integer~$d$ and let
\begin{equation}
\label{eq:p-sos}
g(x) = \sum_{\pi \in \btd} \langle V^{\pi, d}(x),
R^\pi\rangle
\end{equation}
be the polynomial that defines~$\widehat{f}$, where~$R^\pi$ are
positive semidefinite matrices of the appropriate sizes. Notice that
this, together with the construction of the~$V^{\pi, d}$ matrices,
implies that~$g$ is a sum of squares polynomial of degree at most~$2d$
invariant under~$\bt$, and that vice versa all sum of squares
$\bt$-invariant polynomials of degree at most~$2d$ are of this form.

Function~$f$ is the Fourier inverse of~$\widehat{f}$. In
Section~\ref{sec:fourier-transform}, we have seen how the inverse can
be computed when~$\widehat{f}$ is given by an invariant polynomial as
in~\eqref{eq:fhat-definition}. In fact, there is a linear
transformation~$\Fcal\colon \R[x]^{\bt} \to \R[x]^{\bt}$ such that
\[
f(x) = \Fcal[g](x) e^{-\pi\|x\|^2}.
\]
In particular, if~$g$ is given as in~\eqref{eq:p-sos}, then
\[
f(x) = e^{-\pi\|x\|^2} \sum_{\pi \in \btd} \langle \Fcal[V^{\pi,
  d}](x), R^\pi\rangle,
\]
where by applying~$\Fcal$ to a matrix we apply it to each entry and
get a matrix as a result.

With this, we can easily see how to express condition~(i) of
Theorem~\ref{thm:cohn-elkies}. It becomes
\[
\sum_{\pi \in \btd} \langle V^{\pi, d}(0), R^\pi\rangle \geq 1.
\]
The bound provided by the theorem is then~$\vol \Kcal$ times
\[
f(0) = \sum_{\pi \in \btd} \langle \Fcal[V^{\pi, d}](0), R^\pi\rangle;
\]
this will be the objective function of our semidefinite program.

\subsection{Nonpositivity constraint}

We impose the condition~$f(x) \leq 0$ when
$x \notin \Kcal^\circ - \Kcal^\circ$ in two steps by breaking the
domain in which the function has to be nonpositive in two parts, an
unbounded and a bounded one. We then deal with the unbounded part with
an SOS constraint and with the bounded part via sampling.

Let us consider first the unbounded part of the domain. Let~$s$ be a
$\bt$-invariant polynomial such that
\[
\Kcal^\circ - \Kcal^\circ \subseteq \{\, x \in \R^3 : s(x) < 0\,\},
\]
where the set on the right-hand side is bounded.  For instance,
if~$\delta$ is the maximum norm of a vector in~$\Kcal - \Kcal$, then
we may take~$s(x) = \|x\|^2 - \delta^2$.

If there are SOS polynomials~$q_1$ and~$q_2$ such that
\begin{equation}
\label{eq:sos-constraint}
\Fcal[g](x) = -s(x) q_1(x) - q_2(x),
\end{equation}
then~$f$ will be nonpositive in~$\{\, x \in \R^3 : s(x) \geq
0\,\}$.
Now,~$g$ is invariant and hence~$\Fcal[g]$ is invariant. Since~$s$ is
also invariant, we may take both~$q_1$ and~$q_2$ invariant without
loss of generality. So we may use for~$q_1$ and~$q_2$ a
parametrization similar to the one we used for~$g$, but here it is
important to be careful with the choice of degrees of~$q_1$ and~$q_2$.

In principle, the degrees of~$q_1$ and~$q_2$ can be anything as long
as they are high enough so that the identity above may hold. In
practice, it is a good idea to limit the degrees of~$q_1$ and~$q_2$ as
much as possible. For instance, if~$q_2$ is allowed to have a larger
degree than~$\Fcal[g]$, then it is certainly not possible to represent
it in our parametrization with positive \textit{definite} matrices,
and this will make it very difficult to rigorously prove that the
numbers we obtain are indeed bounds.

We fixed the degree of~$g$ to be at most~$2d$ for some
odd~$d$. Then~$\Fcal[g]$ also has degree at most~$2d$. Since~$s$ is
invariant, it has an even degree, say~$2d_s$. We will
impose~$\deg q_1 \leq 2d - 2d_s$ and~$\deg q_2 \leq 2d$. Now we may
parametrize~$q_1$ using positive semidefinite matrices~$S_1^\pi$
for~$\pi \in \btd$ and~$q_2$ using positive semidefinite
matrices~$S_2^\pi$ for~$\pi \in \btd$,
rewriting~\eqref{eq:sos-constraint} as
\[
\sum_{\pi \in \btd} \langle \Fcal[V^{\pi, d}](x),
R^\pi\rangle 
+ \sum_{\pi \in \btd} \langle s(x)V^{\pi, d -
  d_s}(x), S_1^\pi\rangle
+ \sum_{\pi \in \btd} \langle V^{\pi, d}(x),
S_2^\pi\rangle = 0.
l\]
Notice this is a polynomial identity, which should be translated into
linear constraints in our semidefinite program. To do so we need to
express all polynomials in a common basis. A natural choice here is a
basis of the invariant ring~$\R[x]^\bt$, since we work exclusively
with invariant polynomials.

Now we still need to ensure that~$f$ is nonpositive in the bounded set
\[
\Dcal = \{\, x \in \R^3 : s(x) < 0\,\} \setminus (\Kcal^\circ -
\Kcal^\circ).
\]
We do so by using a finite sample of points in~$\Dcal$ and adding
linear constraints requiring that~$\Fcal[g]$ be nonpositive for each
point in the sample. The idea is that, if we select enough points,
then these constraints should ensure that~$f$ is nonpositive
everywhere in~$\Dcal$.

So we choose a finite set~$\Scal \subseteq \Dcal$. Because of the
invariance of~$g$, and hence of~$\Fcal[g]$, we may restrict ourselves
to points in the fundamental domain of~$\bt$ or, in other words, we
may restrict ourselves to points~$(x_1, x_2, x_3)$ with~$0\leq x_1\leq x_2\leq
x_3$. Then we add the constraints
\[
\sum_{\pi\in\btd} \langle \Fcal[V^{\pi,d}](x), R^\pi\rangle \leq 0
\qquad\text{for all~$x \in \Scal$}
\]
to our problem.

\subsection{Full formulation}

Here is the semidefinite programming problem we solve. Recall
that~$d$ is an odd positive integer.
\begin{equation}
\label{eq:full-sdp}
\vcenter{\openup5pt\halign{\hfil#\enspace&$#$\hfil\cr
\text{min}&\displaystyle\sum_{\pi \in \btd} \langle \Fcal[V^{\pi, d}](0), R^\pi\rangle\cr
& \text{(a)} \;\; \displaystyle\sum_{\pi \in \btd} \langle V^{\pi, d}(0), R^\pi\rangle
\geq 1,\cr
& \text{(b)} \;\;\displaystyle\sum_{\pi \in \btd} \langle \Fcal[V^{\pi, d}](x),
R^\pi\rangle 
+ \sum_{\pi \in \btd} \langle s(x)V^{\pi, d -
  d_s}(x), S_1^\pi\rangle\qquad\cr
&\displaystyle\hfill+\sum_{\pi \in \btd} \langle V^{\pi, d}(x),
S_2^\pi\rangle = 0,\cr
& \text{(c)} \;\; \displaystyle\sum_{\pi\in\btd} \langle \Fcal[V^{\pi,d}](x), R^\pi\rangle \leq 0
\qquad\text{for all~$x \in \Scal$},\cr
&\text{$R^\pi$, $S_1^\pi$, and~$S_2^\pi$ are positive
  semidefinite}.\cr
}}
\end{equation}

As mentioned before, to express the SOS constraint above, which is in
fact a polynomial identity, we need to express all polynomials
involved in a given common basis. We use for this a basis
of~$\R[x]^\bt$.

\section{Rigorous verification}
\label{sec:verification}

In this section we discuss how the numerical results obtained can be
turned into rigorous bounds.  For the remainder of this section,
$\Kcal \subseteq \R^3$ will be a convex body such that~$\Kcal - \Kcal$ is
$\mathsf{B}_3$-invariant.

\subsection{Solving the problem and checking the SOS constraint}
\label{sec:solving}

We input problem~\eqref{eq:full-sdp} to a semidefinite programming
solver. In doing so, we are using floating-point numbers to
represent the input data. Solvers also use floating-point numbers
in their numerical calculations, so the solution obtained at the end
is likely not feasible. \emph{But}, if it is close enough to being feasible,
then it can be turned into a feasible solution. To this end it is
important to find a solution in which the minimum eigenvalue of any
matrix~$R^\pi$, $S_1^\pi$ and~$S_2^\pi$ is much larger than the maximum
violation of any constraint.

Here is where it becomes important to make the formulation as tight as
possible, for instance by picking the degrees of polynomials~$g$,
$q_1$, and~$q_2$ in~\eqref{eq:sos-constraint} correctly, so
that~\eqref{eq:full-sdp} admits a \textit{strictly feasible} solution,
that is, a solution in which every matrix is positive definite. It is
also for this reason that we have chosen~$d$ odd, since for even~$d$
the resulting problem is not strictly feasible.

To obtain such a solution we use a two-step approach. First we solve
our problem to get an estimate~$z^*$ on the optimal value. Many
interior point solvers work exclusively with positive definite
solutions, but at the end round the solution to a face of the positive
semidefinite cone. So the resulting solution matrices might have zero
eigenvalues. To overcome this problem, we then pick some small
error~$\eta$ (we usually pick something like~$\eta = 10^{-5}$) and
remove the objective function of the problem, adding it as a
constraint like
\[
\sum_{\pi \in \btd} \langle \Fcal[V^{\pi, d}](0), R^\pi\rangle \leq
z^* + \eta.
\]
So we sacrifice a bit of the optimal value. Most solvers, however,
when dealing with feasibility problems, i.e., problems without an
objective function, return a solution in the \textit{analytic center}
if a solution is found, and that solution will have positive definite
matrices with large minimum eigenvalues. Of course, how large the
minimum eigenvalues will be depends on the choice of~$\eta$.

It is also important to use a solver able to work with high-precision
floating-point numbers. Solvers working with double-precision
floating-point arithmetic have failed to find feasible solutions of
our problem because of numerical stability issues. Moreover, by using
high-precision arithmetic we will get in the end a solution that is
only slightly violated, which is our goal. To solve our problems, we
used the SDPA-GMP solver~\cite{SDPA}.

Say then we have a solution~$(R^\pi, S_1^\pi, S_2^\pi)$ with the
desired property, that is, a solution in which the minimum eigenvalues
are much larger than the maximum constraint violation. Our next step
is to get a bound on the minimum eigenvalue of each matrix
involved. We do it as follows. For each matrix~$A$ in the solution, we
use binary search to find~$\lambda_A > 0$ close to the minimum
eigenvalue of~$A$ so that~$A - \lambda_A I$ has a Cholesky
decomposition~$L L^\tp$. This we do with high-precision floating-point
arithmetic. Then we use instead of~$A$ the
matrix~$\tilde{A} = L L^\tp + \lambda_A I$. We have then a positive
definite matrix and a bound on its minimum eigenvalue. We use interval
arithmetic with high-precision floating-point arithmetic~\cite{Revol2005a}
to represent the new
solution~$(\tilde{R}^\pi, \tilde{S}_1^\pi, \tilde{S}_2^\pi)$ obtained
in this way.

Now we can easily compute how violated the normalization
constraint~(a) of~\eqref{eq:full-sdp} is using interval arithmetic; if
it is violated then we can multiply the solution by a positive number
so as to have it satisfied. It is also easy to compute the objective
value of the solution. We can also use interval arithmetic to compute
an upper bound on the maximum violation of the SOS constraint~(b)
in~\eqref{eq:full-sdp}. To do so, we compute the absolute value of the
coefficient of
\[
r(x) = \sum_{\pi \in \btd} \langle \Fcal[V^{\pi, d}](x),
\tilde{R}^\pi\rangle + \sum_{\pi \in \btd} \langle s(x)V^{\pi, d -
  d_s}(x), \tilde{S}_1^\pi\rangle + \sum_{\pi \in \btd} \langle
V^{\pi, d}(x), \tilde{S}_2^\pi\rangle
\]
with largest absolute value. Here we should note that
matrices~$V^{\pi, m}$ can be expressed using only rationals (if we
work with a basis of~$\R[x]^\bt$ whose elements have only rational
coefficients, as we actually do), and these rationals can be
approximated with interval arithmetic. For the
matrices~$\Fcal[V^{\pi, d}]$ of Fourier inverses we also need
irrationals (namely, powers of~$\pi$), but these can be approximated
with interval arithmetic.

We want to change~$(\tilde{R}^\pi, \tilde{S}_1^\pi, \tilde{S}_2^\pi)$
in order to make~$r$ identically zero. Notice~$r$ is invariant and has
degree up to~$2d$. By construction of the~$V^{\pi,d}$ matrices,~$r$
can be expressed as a linear combination of their entries. In other
words, there are matrices~$T^\pi$ such that
\[
r(x) = \sum_{\pi \in \btd} \langle
V^{\pi, d}(x), T^\pi\rangle.
\]
Then~$(\tilde{R}^\pi, \tilde{S}_1^\pi, \tilde{S}_2^\pi - T^\pi)$
satisfies the SOS constraint~(b). If the numbers in~$T^\pi$ are small
enough compared to the minimum eigenvalue of~$\tilde{S}_2^\pi$,
then~$\tilde{S}_2^\pi - T^\pi$ will be positive semidefinite, and we
will have obtained a solution satisfying the SOS constraint
in~\eqref{eq:full-sdp}. Namely, it suffices to require
\[
\|T^\pi\| \leq \lambda_{\tilde{S}_2^\pi}
\]
for all~$\pi \in \btd$, where~$\lambda_{\tilde{S}_2^\pi}$ is any lower
bound on the minimum eigenvalue of~$\tilde{S}_2^\pi$, which may be
obtained as explained above. Here,
$\|A\| = \langle A, A \rangle^{1/2}$ is the Frobenius norm of
matrix~$A$.

To estimate~$\|T^\pi\|$ we use the following approach in which we do
not explicitly determine $T^\pi$. We find a maximal linearly
independent subset~$\Bcal$ of polynomials inside the set of all
entries of the~$V^{\pi,d}$ matrices for~$\pi \in \btd$.  Now we
create the matrix~$A$ with rows indexed by all monomials occurring in
a polynomial in~$\Bcal$ and columns indexed by~$\Bcal$. An
entry~$(m, a)$ of~$A$, where~$m$ is a monomial and~$a \in \Bcal$,
contains the coefficient of monomial~$m$ in polynomial~$a$. Then we
find a submatrix~$\hat{A}$ of~$A$ consisting of~$|\Bcal|$
linearly-independent rows of~$A$ and we compute~$\hat{A}^{-1}$ using
rational arithmetic. We may compute
\[
\|\hat{A}^{-1}\|_\infty = \max_{i=1,\ldots,|\Bcal|}
\sum_{j=1}^{|\Bcal|} |\hat{A}^{-1}_{ij}|
\]
and observe that the maximum absolute value of any coefficient of the
expansion of~$r$ in basis~$\Bcal$ is at most~$\|\hat{A}^{-1}\|_\infty
\|r\|_\infty$, where~$\|r\|_\infty$ is the maximum absolute value of
any coefficient of~$r$. In this way we may get an estimate
on~$\|T^\pi\|$.

\subsection{Checking the sample constraints}

Even if one uses a great number of sample constraints in condition~(c)
it is unlikely that the resulting function will be nonpositive
in~$\Dcal$ as required. The sample constraints cannot accurately
detect the boundary of~$\Dcal$. However, for some small
factor~$\alpha > 1$, which we hope will be small if the sample was
fine enough, the function will be nonpositive in the domain
\[
\Dcal' = \{\, x \in \R^3 : s(x) < 0\,\} \setminus \alpha (\Kcal^\circ -
\Kcal^\circ).
\]

One may quickly estimate a good value for~$\alpha$ by testing the
function on a fine grid of points. Then all that is left to do is
check that the function is really nonpositive in~$\Dcal'$. Of course,
the larger the~$\alpha$, the worse the bound will be because it needs
to be multiplied by~$\alpha^3$.

\begin{figure}[htb]
\includegraphics[width=5.5cm]{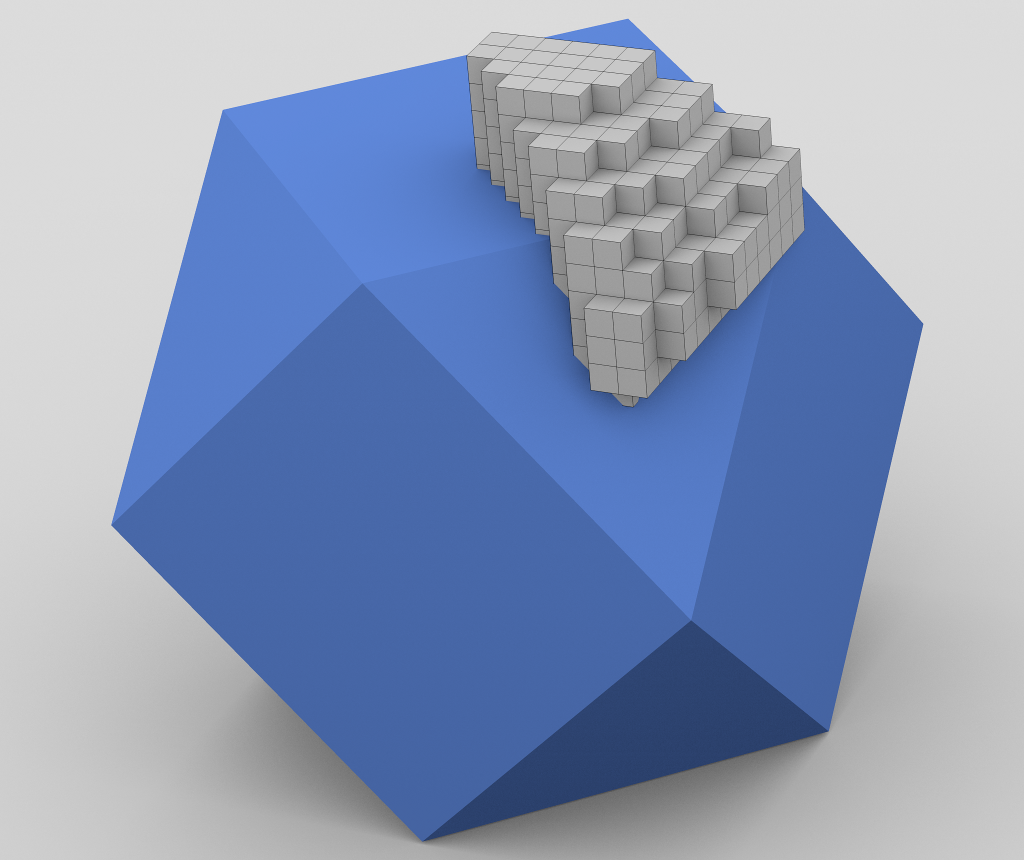}
\hskip1cm
\includegraphics[width=5.5cm]{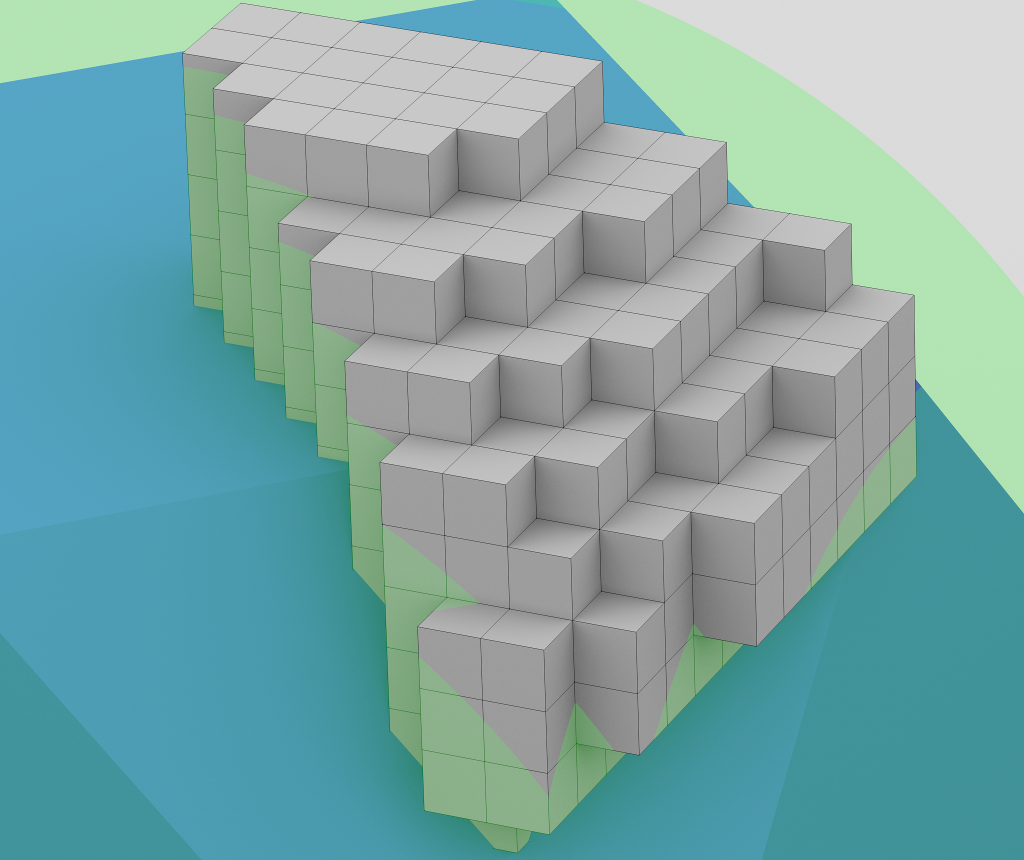}

\caption{Here we have an initial cube partition for~$\Dcal'$
  where~$\Kcal$ is a regular tetrahedron so that~$\Kcal - \Kcal$ is
  the cuboctahedron with circumradius~$1$. We take~$\alpha = 1.02$
  and~$s(x) = \|x\|^2 - 1$. On the left we show the whole
  cuboctahedron and the partition. On the right we show the partition
  in detail; the unit sphere is also shown in green.}
\label{fig:partition}
\end{figure}

In our approach we use interval arithmetic to evaluate the
polynomial~$\Fcal[g]$, so as to obtain rigorous results. We consider a
partition of~$\R^3$ into cubes of side-length~$\delta$ for some
small~$\delta$ and we let~$\Ccal$ be the set of all partition cubes
that contain at least one point~$(x_1, x_2, x_3) \in \Dcal'$
with~$0\leq x_1 \leq x_2 \leq x_3$. Note that~$\Ccal$ is finite and covers~$\Dcal'$;
Figure~\ref{fig:partition} shows an example initial partition
when~$\Kcal$ is the regular tetrahedron.

We then check that for every point in~$\bigcup_{C \in \Ccal} C$ the
polynomial~$\Fcal[g]$ is nonpositive. We do that as follows.

First, for every cube~$C \in \Ccal$ we compute an upper bound of the
norm of the gradient of $\Fcal[g]$, a number~$\nu_C$ such that
\begin{equation}
\label{eq:gradientupperbound}
\|\nabla \Fcal[g](x)\| \leq \nu_C
\end{equation}
for all~$x \in C$. This is easy to do with interval arithmetic. We
have the coefficients of~$\Fcal[g]$ represented by intervals. A
cube~$C = [x_1, y_1] \times [x_2, y_2] \times [x_3, y_3]$ is the
product of three intervals. We then only have to
compute~$\nabla \Fcal[g]([x_1, y_1], [x_2, y_2], [x_3, y_3])$ using
interval arithmetic. This will give us a vector of intervals~$([l_1,
u_1], [l_2, u_2], [l_3, u_3])$ such that for all~$(x_1, x_2, x_3) \in
C$ we have
\[
(l_1, l_2, l_3) \leq \nabla \Fcal[g](x_1, x_2, x_3) \leq (u_1, u_2, u_3).
\]
From this it is easy to compute a number~$\nu_C$
satisfying~\eqref{eq:gradientupperbound}.

Next, for a fixed integer~$N \geq 1$, say, we uniformly divide each
side of the cube~$C$ into~$N$ intervals, obtaining a grid of points
inside of~$C$. In other words, if~$x_C$ is the lower-left corner
of~$C$, we consider the set of points
\[
C_N = \{\, x_C + (a, b, c) \delta / N : 0 \leq a, b, c \leq N,\, a,b,c
\in \N\,\}.
\]
At least one point of~$C$ belongs to~$\Dcal'$, and hence at least one
point of~$C_N$ is not in~$\alpha(\Kcal^\circ - \Kcal^\circ)$. Let
then~$d(C,N)$ be the maximum minimum distance from any point
of~$C \setminus \alpha(\Kcal^\circ - \Kcal^\circ)$ to a point
of~$C_N \setminus \alpha(\Kcal^\circ - \Kcal^\circ)$ and let
\[
\mu(C, N) = \max\{\, \Fcal[g](x) : x \in C_N \setminus
\alpha(\Kcal^\circ - \Kcal^\circ)\,\}.
\]
If~$\mu(C, N) > 0$, then the function is not nonpositive in the
required domain. We hope however that, for our choice of~$\alpha > 1$,
we will have~$\mu(C, N) < 0$. Suppose that this is the case.
Given~$x \in C \setminus \alpha(\Kcal^\circ - \Kcal^\circ)$, let~$x'$
be the point in~$C_N \setminus \alpha(\Kcal^\circ - \Kcal^\circ)$
closest to~$x$. By the mean-value theorem we have that
\[
|\Fcal[g](x) - \Fcal[g](x')| \leq \nu_C \|x - x'\| \leq \nu_C d(C, N).
\]
So, if
\begin{equation}
\label{eq:suff}
\nu_C d(C, N) \leq |\mu(C, N)| \leq |\Fcal[g](x')|,
\end{equation}
then~$\Fcal[g](x) \leq 0$. Checking condition \eqref{eq:suff} for all
$x' \in C_N \setminus \alpha(\Kcal^\circ - \Kcal^\circ)$ gives us a sufficient condition that
allows us to conclude that~$\Fcal[g](x) \leq 0$ for
all~$x \in C \setminus \alpha(\Kcal^\circ - \Kcal^\circ)$.

We still have to estimate~$d(C, N)$, but that is a simple
matter. There are two cases.
If~$C \cap \alpha(\Kcal^\circ - \Kcal^\circ) \neq \emptyset$,
then~$d(C, N) \leq (\delta / N)\sqrt{3}$; if not,
then~$d(C, N) \leq (\delta/(2 N)) \sqrt{3}$.

So our strategy is to process each cube in~$\Ccal$. For each cube, we
start with~$N = 2$ and check if that is enough to conclude that the
function is nonpositive in the cube. If not, then we
increase~$N$. Once all cubes have been processed, we know that the
function is nonpositive everywhere in the domain. Finally, notice that
we always use interval arithmetic to perform all computations, thus
obtaining rigorous results at the end, once the procedure terminates.

There is only one extra issue that is conceptually simple but that
makes things technically harder. Computing with interval arithmetic is
very slow, and hence if too many cubes would require dense grids (say,
with hundreds of points per side), the computation would take several
months. The size of the grid required by a cube is however directly
proportional to the upper bound on the norm of the gradient, which is
better the smaller the cube is. So, by taking smaller~$\delta$, we can
improve on the grid sizes. But by changing~$\delta$ globally, we
increase the total number of cubes, possibly slowing down the total
computation time.

A better strategy is as follows: if the grid size required by a cube
is greater than a certain threshold (we use~$30$), then we split the
cube at its center creating eight new cubes and keeping only those
that intersect the domain. Then we process the resulting cubes
instead, which are smaller and therefore lead to better grid
sizes. This splitting process is carried out recursively, up to a
certain maximum depth.

Finally, when one estimates the required grid size it may happen that,
from one iteration to another, the grid size~$N$ is increased only
slightly. This should be avoided, since computing the function is
quite expensive. So our approach is as follows: first, we carry out
the whole verification procedure using double-precision floating-point
arithmetic for function evaluation, but not for the other
computations. This is quite fast, finishing in a few hours. Then we
use the estimated grid size for each cube in a checking routine that
remakes all calculations using interval arithmetic.

\subsection{Further implementation details}

Table~\ref{tab:verified} contains the list of bounds we computed and
rigorously verified using the approach described in this section. 

The solutions, as well as the verification scripts and programs, can be found as ancillary files from the \texttt{arXiv.org} e-print archive.

\begin{table}[htb]
\begin{center}
\begin{tabular}{lcc}
Body&Upper bound& Factor $\alpha$\\
\hline
Regular octahedron ($B_3^1$)       & $0.972912750$   & $1.001$  \\
$B_3^3$                      & $0.823611150$   & $1.002$  \\
$B_3^4$                           & $0.874257405$   & $1$      \\
$B_3^5$                           & $0.922441815$   & $1.005$  \\
$B_3^6$                            & $0.933843309$   & $1$      \\
Regular tetrahedron              & $0.374568355$   & $1.02$   \\
Truncated cube                    & $0.984519783$   & $1.003$ \\
Truncated tetrahedron          & $0.729209804$   & $1.023$ \\
%Rhombic cuboctahedron             & $0.879465169$   & $1.009$ \\
%Truncated cuboctahedron           & $0.884572870$   & $1.018$\\
\hline
\end{tabular}
\end{center}
\bigskip
\caption{List of rigorous bounds together with the factor
  $\alpha$ we needed in the verification.}
\label{tab:verified}
\end{table}

The procedure to verify the SOS constraints was implemented as a
Sage~\cite{Sage} script \texttt{verify.sage} and runs in Sage~6.2; see
the documentation file \texttt{README\_SOSChecking}. The approach to test
that the function is nonpositive in the domain was implemented as a
C++11 program called {\tt checker} using the MPFI
library~\cite{Revol2005a} for interval arithmetic. Verification time
was in all but one case under~$2$ days; in the case of the regular
octahedron it took several weeks. More documentation of the C++11
program can be found in \texttt{README\_SampleChecking} and a description
of the classes in \texttt{docu.pdf}.

\begin{figure}[htb]
\includegraphics[width=12cm]{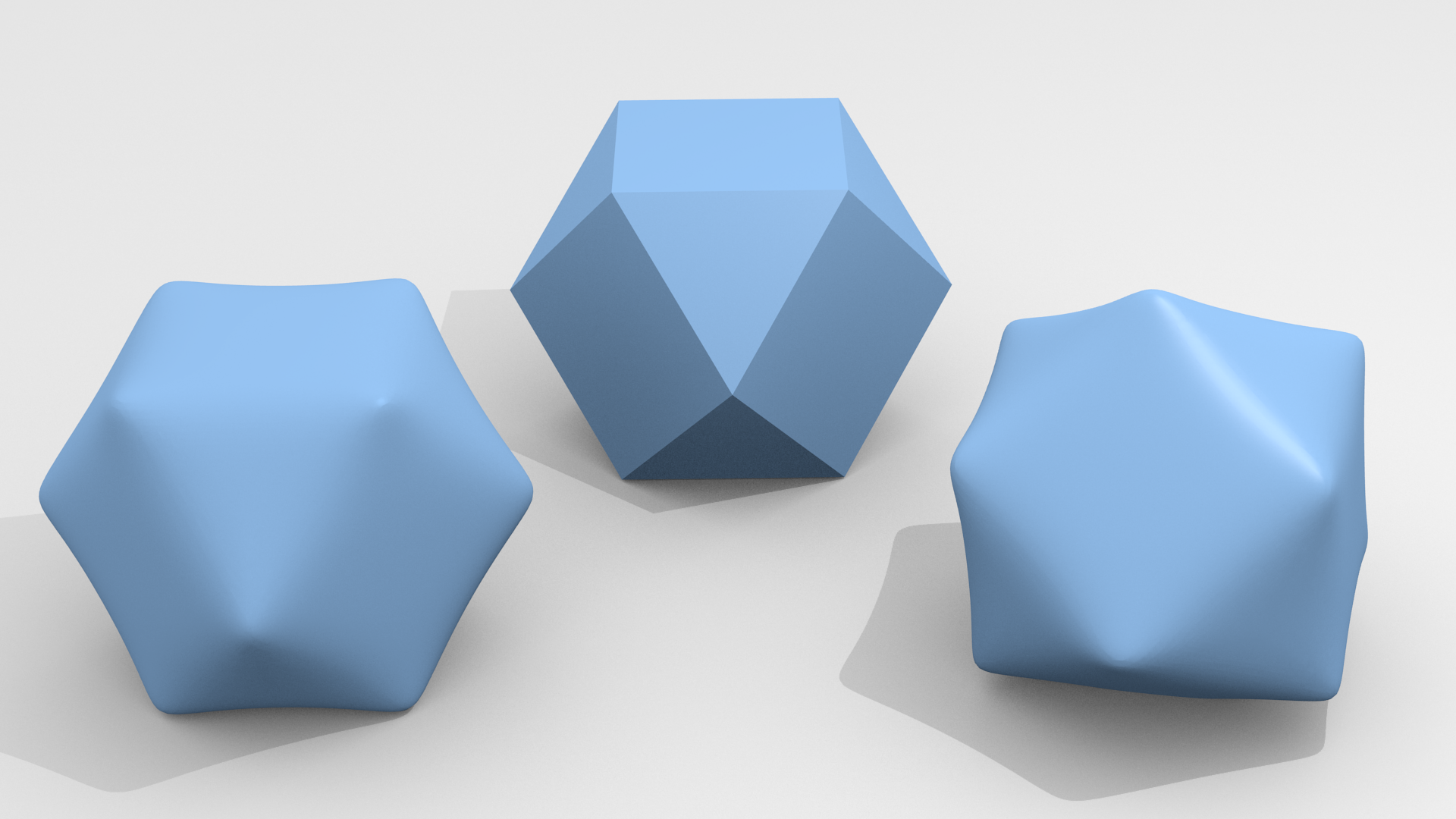}

\caption{The cuboctahedron (center) with two rotations of the set
  $\{\,x \in \R^3 : \Fcal[g](x) \geq 0\,\}$, where~$g$ is the
  polynomial given by the solution to our problem for the regular
  tetrahedron.}
\label{fig:cuboc}
\end{figure}

It is interesting to observe that the polynomials~$\Fcal[g]$ obtained
from the solutions to the semidefinite programs we consider provide
interesting low degree polynomial approximations of the Minkowski
difference $\Kcal - \Kcal$; We used $d = 13$ in our computations,
so that $g$ and $\mathcal{F}[g]$ have degree $26$. Figure~\ref{fig:cuboc} shows the
cuboctahedron, which is the Minkowski difference of two regular
tetrahedra, and the
region~$\{\, x \in \R^3 : \Fcal[g](x) \geq 0\, \}$, where~$g$ is given
by the solution to our problem for the tetrahedron. Notice how the
region approximates the cuboctahedron. In fact, the upper bounds we
computed are also bounds for translative packings of the nonconvex
bodies determined by these polynomial approximations.

\section*{Acknowledgements}

The fourth author thanks Peter Littelmann for a helpful discussion. We
also thank the referees for their thorough comments which helped to
improve the paper.

\end{document}